\documentclass [a4paper,twoside,13pts]{article}
\usepackage [francais]{babel}
\usepackage [latin1] {inputenc}
\usepackage [T1] {fontenc}
\usepackage{amssymb}
\usepackage{amsmath}
\usepackage{stmaryrd}
\usepackage{graphicx}
\date {18th January 2007}
\title{Penalizations of the Brownian motion by a functional of its local times }
\author {Joseph Najnudel}
\begin{document} 
\maketitle 
\noindent
\textbf{Abstract : } In this article, we study the family of
probability measures (indexed by $t \in \mathbf{R}_+^*$), obtained by
penalization of the Brownian motion by
a given functional of its local times at time $t$. \\  We prove that
this family tends to a limit measure when $t$ goes to infinity if
the functional satisfies some conditions of domination, and we check
these conditions in several particular cases.  
\section*{Introduction}
\noindent
Brownian penalizations have been studied in several articles, in
particular
in \cite{6}, \cite{7}, \cite{8}. The general
principle of these penalizations is the following : let $\mathbf{W}$
be the Wiener measure on $\mathcal{C} (\mathbf{R}_+, \mathbf{R})$,
$(X_t)_{t \geq 0}$ the canonical process, and $(\Gamma_t)_{t \geq 0}$
a family of positive weights such that $0 < \mathbf{W} [ \Gamma_t] <
\infty$; we consider the family of probability measures 
$(\mathbf{W}_t)_{t \geq 0}$, obtained from $\mathbf{W}$, by
``penalization'' with the weight $\Gamma$ : 
$$\mathbf{W}_t = \frac{\Gamma_t}{\mathbf{W}
  [\Gamma_t]} . \mathbf{W}$$
\noindent
In many different particular cases, the family $(\mathbf{W}_t)_{t \geq
  0}$ tends to a limit measure $\mathbf{W}_{\infty}$ as $t \rightarrow
  \infty$, in the following
  sense : for all $s \geq 0$, and for $\Lambda_s$ measurable with
  respect to $\mathcal{F}_s = \sigma \{X_u, u \leq s \}$ : $$ \mathbf{W}_t (\Lambda_s)
  \underset{t \rightarrow \infty}{\rightarrow} \mathbf{W}_{\infty} (\Lambda_s)$$ 
\noindent 
Up to now, there does not exist a general theorem which covers all the
different cases for which convergence holds. On the other hand, we remark that in
many of 
these cases, one has : $$ \Gamma_t = F((l_t^y(X))_{y \in
  \mathbf{R}})$$ where $(l_t^y(X))_{y \in
  \mathbf{R}}$ is the family of the local times of $(X_s)_{s \leq t}$,
and $F$ is a measurable functional from $\mathcal{C} (\mathbf{R},
\mathbf{R}_+)$ to $\mathbf{R}_+$. \\ \\
These two facts led us to prove that if $\Gamma$ is of this form, the
limit measure $\mathbf{W}_{\infty}$ exists for a ``large'' class of
functionals $F$. \\ \\
This proof is the main topic of our article, which is divided into
six sections. \\ \\ In the first one, we define
and explain the  notations we need to prove our main theorem,
which is stated at the end of the section. \\ \\
In Section 2, we prove an equality satisfied by an approximation of
a given functional of local times, and in Section 3, we majorize the error term
corresponding to this approximation. \\ \\
This allows us to obtain, in Section 4, the asymptotic behaviour of the expectation of
functionals which satisfy some particular conditions, and finally we prove the
main theorem in Section 5. \\ \\
In Section 6, we study the four following examples, for which the
Theorem applies : \\ \\
1) $F((l^y)_{y \in \mathbf{R}}) = \phi(l^0)$ (which corresponds to
$\Gamma_t = \phi(l_t^0(X))$), where $\phi$ is a function from
$\mathbf{R}_+$ to $\mathbf{R}_+$, dominated by an integrable and 
decreasing function $\psi$. \\ \\
2) $ F((l^y)_{y \in \mathbf{R}}) = \phi ( \inf \{ y \geq 0, l^y = 0
\} )$ (which corresponds to the weight \\  $\Gamma_t = \phi( \sup
\{ X_s, s \leq t \})$), where $\phi$ is a function from $\mathbf{R}_+
\cup \{ \infty \}$ to $\mathbf{R}_+$, dominated by a decreasing
function $\psi$, which is integrable on $\mathbf{R}$. \\ \\
3)  $ F((l^y)_{y \in \mathbf{R}}) = \exp \left( -
  \int_{-\infty}^{\infty} V(y) l^y \, dy \right)$, where $V$ is a
positive measurable function, not a.e. equal to zero, and integrable
with respect to $(1+y^2) dy$. \\ \\
4)  $ F((l^y)_{y \in \mathbf{R}}) = \phi(l^{y_1}, l^{y_2})$, where
$y_1<y_2$ and $\phi(l_1,l_2) \leq h(l_1 \wedge l_2)$, for a decreasing
and integrable function $h$. \\ \\
The three first examples have been already studied by B. Roynette,
P. Vallois and M. Yor. \\ \\
As a help to the reader, we mention that Sections 2 and 3 are quite
technical, but it is
possible to read the
details of these sections after Sections 4 and 5, which contain the
principal steps of the proof of the Theorem.    
\section{Notations and statement of the main theorem}
\noindent
In this article,  $(B_t)_{t \geq 0}$ denotes a standard
one-dimensional Brownian motion,  $(L_t^y)_{t \geq 0, y \in \mathbf{R}}$
the bicontinuous version of its local times, and $(\tau_l^a)_{l \geq 0, a
  \in \mathbf{R}}$ the family of its inverse local times. \\ \\
To simplify these notations, we put  $T_a = \tau_0^a$ (first hitting
time at $a$ of $B$) and  $\tau_l^0=\tau_l$.  \\ \\
For every $l \in \mathbf{R}_+$, $(Y_{l,+}^y)_{y \in \mathbf{R}}$ denotes a random
process defined on the whole real line, such that its ``positive part'' $(Y_{l,+}^y)_{y \geq 0}$ is a 2-dimensional squared
Bessel process (BESQ(2)), its ``negative part'' $(Y_{l,+}^{-y})_{y \geq 0}$ is an
independent 0-dimensional squared Bessel process (BESQ(0)), and
its value at zero $Y_{l,+}^0$ is equal to $l$. In particular, by
classical properties of BESQ(0) and BESQ(2) processes, there exists a.s. $y_0 \leq 0$
such that $Y_{l,+}^y = 0$  iff $y \leq y_0$. \\
We define also  $(Y_{l,-}^y)_{y \in \mathbf{R}}$ as a
process which has the same law as $ (Y_{l,+}^{-y})_{y \in
  \mathbf{R}}$, the process obtained from $ (Y_{l,+}^{y})_{y \in
  \mathbf{R}}$ by ``reversing the time''. \\ \\
In one of the penalization results shown in \cite{8}, B. Roynette,
P. Vallois and M. Yor obtain a limit process
$(Z_t^l)_{t \geq 0}$, such that $Z_t^l = B_t$ for $t \leq \tau_l$,
$(|Z_{\tau_l + u}^l|)_{u \geq 0}$ is a BES(3) process independent
of $B$, and $\epsilon = \operatorname{sgn} (Z_{\tau_l + u}^l)$ ($u > 0$) is an
independent variable such that $\mathbf{P} (\epsilon = 1) = \mathbf{P}
(\epsilon = -1) = 1/2 $. This process can be informally considered to
be a Brownian motion conditionned to have a total local time equal to
$l$ at level zero. By applying Ray-Knight theorems for Brownian local
times (see \cite{5}) to $(Z_t^l)_{t \geq 0}$, it
is possible to show that the law of the family of its total local
times is the half-sum of the laws of $ (Y_{l,+}^{y})_{y \in
  \mathbf{R}}$ and $ (Y_{l,-}^{y})_{y \in
  \mathbf{R}}$ ($ (Y_{l,+}^{y})_{y \in
  \mathbf{R}}$ corresponds to the paths of $(Z_t^l)_{t \geq 0}$ such
that $\epsilon = 1$, and  $ (Y_{l,-}^{y})_{y \in
  \mathbf{R}}$ corresponds to the paths such that $\epsilon = -1$). \\
This explains why the processes  $ (Y_{l,+}^{y})_{y \in
  \mathbf{R}}$ and  $ (Y_{l,-}^{y})_{y \in
  \mathbf{R}}$ occur naturally in 
the description of the asymptotic behaviour of Brownian local times. \\ \\
We also need to define some modifications of $(Y_{l,+}^y)_{y \in
  \mathbf{R}}$  and  $(Y_{l,-}^y)_{y \in \mathbf{R}}$ : for $l \geq
0$, $a \geq 0$,  $(Y_{l,a}^y)_{y \in \mathbf{R}}$ denotes a
process such that $(Y_{l,a}^y)_{y \geq 0}$ is markovian with the
infinitesimal generator of BESQ(2) for $y \leq
a$ and the infinitesimal generator of BESQ(0) for $y \geq a$,
$(Y_{l,a}^{-y})_{y \geq 0}$  is an independent BESQ(0) process, and
$Y_{l,a}^0 = l$. For $a \leq 0$, $(Y_{l,a}^y)_{y \in \mathbf{R}}$ has
the same law as $(Y_{l,-a}^{-y})_{y \in \mathbf{R}}$. \\ \\ 
Now, let $F$ be a functional from $\mathcal{C} (\mathbf{R},
\mathbf{R}_+)$ to $\mathbf{R}_+$, which is measurable with respect to the
$\sigma$-field generated by the topology of uniform convergence on compact
sets. We
consider the following quantities, which will naturally appear in the 
asymptotics of $\mathbf{E} [F((L_t^y)_{y \in \mathbf{R}})]$ :

$$ I_+ (F) = \int_0^{\infty} dl \, \mathbf{E} [F((Y_{l,+}^y)_{y \in
  \mathbf{R}})] $$
$$ I_- (F) = \int_0^{\infty} dl \, \mathbf{E} [F((Y_{l,-}^y)_{y \in
  \mathbf{R}})] $$
$$I(F)=I_+(F) +  I_-(F)$$
\noindent
We observe that $I(F)$ is the integral of $F$ with respect to the
$\sigma$-finite measure $I$ on $\mathcal{C} (\mathbf{R},
\mathbf{R}_+)$, defined by : $$ I = \int_0^{\infty} dl \,
P_{l,+} + \int_0^{\infty} dl \,
P_{l,-} $$ where $P_{l,+}$ is the law of  $(Y_{l,+}^y)_{y \in
  \mathbf{R}}$ and $P_{l,-}$ is the law of  $(Y_{l,-}^y)_{y \in
  \mathbf{R}}$. \\ \\
At the end of this section,  we give some conditions on $F$ which turn out to be sufficient to
obtain our penalization result. \\
Unfortunately, these conditions are not very simple and we
need three more definitions before stating the main Theorem : \\ \\
\textbf{Definition 1 } (a condition of domination) :  \textit{Let $c$ and $n$ be in $\mathbf{R}_+$
(generally $n$ will be an integer). For every decreasing
function $h$ from $\mathbf{R}_+$ to $\mathbf{R}_+$, we say that a
measurable function $F$ from $\mathcal{C}
(\mathbf{R}, \mathbf{R}_+)$ to $\mathbf{R}_+$ satisfies the condition
$C(c,n,h)$ iff the following holds for every continuous function $l$
from $\mathbf{R}$ to $\mathbf{R}_+$ : \\ \\
1) $F((l^y)_{y \in \mathbf{R}})$ depends only on $(l^y)_{y \in
  [-c,c]}$. \\ \\
2)  $F((l^y)_{y \in \mathbf{R}}) \leq \left( \frac{ \underset{y \in
      [-c,c]}{\sup} l^y + c}{ \underset{y \in
      [-c,c]}{\inf} l^y + c} \right) ^n h \left(  \underset{y \in
      [-c,c]}{\inf} l^y \right) $}  \\ \\ \\
Intuitively, a functional of the local times satisfies the above
condition if it depends only on the local times on a compact set, and
if it is small when these local times are large and
don't vary too much. \\ \\
Now, let us use the notation :  $$N_c(h) = c h(0) + \int_0^{\infty}
h(y) dy$$ \noindent
If $N_c(h) < \infty$, it is possible to prove our main theorem for all 
functionals $F$ which satisfies the condition $C(c,n,h)$, but this
condition is restrictive, since the functional $F$ must not depend on
the local times outside of  $[-c,c]$. \\ \\
In order to relax this restriction, we need the following definition :
\\ \\
\textbf{Definition 2 } (a less restrictive condition of domination) :  
\textit{ Let $n$ be in $\mathbf{R}_+$ and $F$ be a
 positive and measurable function from $\mathcal{C} (\mathbf{R},
 \mathbf{R}_+ )$ to $\mathbf{R}$. \\ 
For all $M \geq 0$, let us say that $F$ satisfies the condition
 $D(n,M)$ iff there exists a sequence $(c_k)_{k \geq 1}$ in $[1,
 \infty[$, a sequence $(h_k)_{k \geq 1}$ of decreasing
 functions from $\mathbf{R}_+$ to $\mathbf{R}_+$, and a sequence
 $(F_k)_{k \geq 0}$ of measurable functions from $\mathcal{C} (
\mathbf{R}_+, \mathbf{R} )$ to $\mathbf{R}_+$, such that : \\ \\
1) $F_0 = 0$ and $(F_k)_{k \geq 1}$ tends to $F$ pointwise. \\ \\
2) For all $k \geq 1$, $|F_k - F_{k-1}|$ satisfies the condition
 $C(c_k,n,h_k)$. \\ \\
3) $\underset{k \geq 1}{\sum} N_{c_k} (h_k) \leq M$. \\ \\ 
We define the quantity $N^{(n)} (F)$ as the infimum of $M \geq 0$ such
 that $F$ satisfies the condition $D(n,M)$.} \\ \\ \\
Intuitively, if $N^{(n)}(F) < \infty$, it means that $F$ can be
 well-approximated by functionals which satisfy conditions given in
 Definition 1. \\ 
In particular,  if $F$ satisfies the condition $C(c, n, h)$ for $c
 \geq 1$, one has : $ N^{(n)} (F) \leq N_c(h)$ (one can prove that $F$
 satisfies the
 condition $ D(n, N_c(h))$, by taking in Definition 2 : $c_k=c$, $h_k =
 h \mathbf{1}_{k=1}$, $F_0 = 0$ and $F_k = F$ if $k \geq 1$). \\ \\ \\
Now, for a given functional $F$, we need to define some other fonctionals, 
informally obtained from $F$  by ``shifting'' the space and adding a given
function to the local time family. \\ \\
More precisely, let us consider the following definition : \\ \\
\textbf{Definition 3 } (local time and space shift) : \textit{Let $x$ be a real number. If $F$ is a
measurable functional from $\mathcal{C}(\mathbf{R}_+, \mathbf{R})$ to
$\mathbf{R}_+$, and if $(l_0^y)_{y \in
  \mathbf{R}}$ is a continuous function from $\mathbf{R}$ to
$\mathbf{R}_+$, we denote by $F^{(l_0^y)_{y \in \mathbf{R}}, x}$ the functional
from $\mathcal{C} (\mathbf{R}, \mathbf{R}_+)$ to $ \mathbf{R}_+$
which satisfies : $$F^{(l_0^y)_{y \in \mathbf{R}}, x} ((l^y)_{y \in
  \mathbf{R}}) = F((l_0^y + l^{y-x})_{y \in \mathbf{R}})$$ for every
function $(l^y)_{y \in \mathbf{R}}$.} \\ \\ \\
This notation and the functionals defined in this way appear naturally
when we consider
the conditional expectation : $\mathbf{E} [F((L_t^y)_{y \in
  \mathbf{R}}) | (B_u)_{u \leq s}]$, for $0<s <t$, and apply the Markov
property. \\ \\ \\
We are now able to state the main theorem of the article : \\ \\
\textbf{Theorem : } \textit{Let $F$ be a functional from $\mathcal{C}
(\mathbf{R}, \mathbf{R}_+)$ to $\mathbf{R}_+$ such that $I(F) > 0$ and
$N^{(n)} (F) < \infty$ for some $n \geq 0$. \\ 
If $\mathbf{W}$ denotes the standard Wiener measure on $\mathcal{C}
(\mathbf{R}_+, \mathbf{R})$, $(X_t)_{t \geq 0}$ the canonical process, and
$(l_t^y (X))_{t \in \mathbf{R}_+, y \in \mathbf{R}}$  the continuous
family of its local times ($\mathbf{W}$-a.s. well-defined), the
probability measure : $$\mathbf{W}_t^F = \frac{F \left( \left( l_t^y
      (X) \right)_{y \in \mathbf{R}} \right)}{ \mathbf{W} \left[ F
    \left( \left( l_t^y (X) \right)_{y \in \mathbf{R}} \right) \right]
} . \mathbf{W}$$ 
is well-defined for every $t$ which is large enough, and there exists a
probability measure  $\mathbf{W}_{\infty}^F$ such that : $$
\mathbf{W}_t^F (\Lambda_s) \underset{t \rightarrow
  \infty}{\rightarrow} \mathbf{W}_{\infty}^F  (\Lambda_s)$$
for every $s \geq 0$ and $\Lambda_s \in \mathcal{F}_s = \sigma \{X_u, u \leq s \}$. \\
\\
Moreover, this limit  measure satisfies the following equality : 
$$  \mathbf{W}_{\infty}^F  (\Lambda_s) =  \mathbf{W} \left(
  \mathbf{1}_{\Lambda_s} . \frac{ I \left(  F^{(l_s^y (X))_{y \in
  \mathbf{R}}, X_s} \right) }{I(F)} \right) $$}
\noindent 
\textbf{Remark 1.1 : } A consequence of the Theorem is the fact that if
$I(F) > 0$ and $N^{(n)} (F) < \infty$ for some $n \geq 0$, the process
$\frac{(I
(F^{(L_s^y)_{y \in \mathbf{R}}, B_s}))_{s \geq 0}}{I(F)}$ is a
martingale. In three of the four examples studied in Section 6, we
compute explicitly this martingale, and in the two first ones, we
check that this computation agrees with the results obtained by B. Roynette,
P. Vallois and M. Yor. \\ \\
\textbf{Remark 2.1 : } We point out that our notation,
$l_t^y(X)$, for the local times given in the Theorem, differs from
the notation $L_t^y$, which is used for the local times
of $(B_s)_{s \leq t}$. This is because, in one
case, we consider the canonical process $(X_t)_{t \geq 0}$ on a given
probability space, and in the other case, we consider a Brownian motion
on a space which is not made precise. Hence, the two mathematical objects
 deserve different writings, despite the fact that
they are strongly related. 
\section{An approximation of the functionals of local times} 
\noindent
In order to prove the Theorem, we need to study the expectation of
$F((L_t^y)_{y \in \mathbf{R}})$, where $F$ is a function from
$\mathcal{C} (\mathbf{R}, \mathbf{R}_+)$ to $\mathbf{R}_+$. \\ \\
However, in general, it is difficult to do that directly, so in this section, we
will replace  $F((L_t^y)_{y \in \mathbf{R}})$ by an approximation.  \\ \\
For the study of this approximation, we need to consider the following
quantities :  
$$ \mathcal{I}_{l,+}^c = \int_{-c}^c Y_{l,+}^y dy, \, \, \, \, 
\mathcal{I}_{l,-}^c = \int_{-c}^c Y_{l,-}^y dy, \, \,  \, \,
\mathcal{I}_{l,a}^c = \int_{-c}^c Y_{l,a}^y dy$$
\noindent
for $c \in \mathbf{R}_+$ or $c = \infty$, $a \in \mathbf{R}$;
$$ \mathcal{Y}_{l,+}^c= \frac{1}{2}( Y_{l,+}^c + Y_{l,+}^{-c}), \, \,
\, \,  
\mathcal{Y}_{l,-}^c= \frac{1}{2}( Y_{l,-}^c + Y_{l,-}^{-c}),   \, \,
\, \,  
\mathcal{Y}_{l,a}^c= \frac{1}{2}( Y_{l,a}^c + Y_{l,a}^{-c})$$
\noindent 
for $c \in \mathbf{R}_+$,  $a \in \mathbf{R}$; 
$$I_{c,t,+} (F) = \int_0^{\infty} dl \, \mathbf{E} \left[ F((Y_{l,+}^y)_{y \in
  \mathbf{R}}) \frac{e^{-(\mathcal{Y}_{l,+}^c)^2/2(t -
  \mathcal{I}_{l,+}^c)}}{\sqrt{1-\mathcal{I}_{l,+}^c/t}} \phi \left(
  \frac{\mathcal{I}_{l,+}^c}{t} \right) \right]$$
$$I_{c,t,-} (F) = \int_0^{\infty} dl \, \mathbf{E} \left[ F((Y_{l,-}^y)_{y \in
  \mathbf{R}}) \frac{e^{-(\mathcal{Y}_{l,-}^c)^2/2(t -
  \mathcal{I}_{l,-}^c)}}{\sqrt{1-\mathcal{I}_{l,-}^c/t}} \phi \left(
  \frac{\mathcal{I}_{l,-}^c}{t} \right) \right]$$ \noindent 
and
$$I_{c,t} (F) = I_{c,t,+} (F) + I_{c,t,-} (F)$$
\noindent
for  $c \in \mathbf{R}_+$, $t > 0$, where  $\phi$ denotes the function
from $\mathbf{R}_+$ to $\mathbf{R}_+$
such that $\phi(x) = 1$ in $x \leq 1/3$, $\phi(x) = 2-3x$ if $1/3 \leq
x \leq 2/3$ and $\phi(x) = 0$ if $x \geq 2/3$ (in particular, this
function is continuous with compact support included in $[0,1[$). \\ \\
We observe that the expression $ \frac{e^{(\mathcal{Y}_{l,+}^c)^2/2(t -
  \mathcal{I}_{l,+}^c)}}{\sqrt{1-\mathcal{I}_{l,+}^c/t}}$ is not
well-defined if $\mathcal{I}_{l,+}^c \geq t$; but this is not important
here, since $\phi(\mathcal{I}_{l,+}^c / t)=0$ in that case. \\ \\
Now, the main result of this section is the following proposition : \\ \\ 
\textbf{Proposition 2 : } \textit{For all measurable functionals from
$\mathcal{C} (\mathbf{R_+}, \mathbf{R})$ to $\mathbf{R}_+$, such that
$F((l^y)_{y \in \mathbf{R}})$ depends only on $(l^y)_{y \in [-c,c]}$
for some $c \geq 0$, the following equality holds : 
$$\sqrt{2 \pi t} \, \mathbf{E} \left[ F((L_t^y)_{y \in \mathbf{R}})
  \mathbf{1}_{|B_t| \geq c} \phi \left( \frac{1}{t} \int_{-c}^c L_t^y
  dy \right) \right] = I_{c,t} (F)$$
for all $t >0$.} \\ \\
\noindent
\textbf{Proof : } Let $G_0$ be a functional from $\mathcal{C}
  (\mathbf{R}_+, \mathbf{R}) \times \mathbf{R}_+$ to $\mathbf{R}_+$,
  such that the process : $(G_0((X_s)_{s \geq 0}, t))_{t \geq 0}$,
  defined on the canonical space $\mathcal{C}
  (\mathbf{R}_+, \mathbf{R})$, is progressively measurable. \\ \\
For every continuous function $\omega$ from $\mathbf{R}_+$ to
  $\mathbf{R}$, $G_0( (\omega_s)_{s \geq 0}, t)$ depends only on
  $(\omega_s)_{s \leq t}$; let us take : $$G( (\omega_s)_{s \leq t}) =
  G_0( (\omega_s)_{s \geq 0}, t)$$
\noindent
Now, by results by C. Leuridan (see \cite{2}), P. Biane
  and M. Yor (see \cite{1}), one has : 
$$\int_0^{\infty} dt \, G((B_s)_{s \leq t}) = \int_0^{\infty} dl \int_{-
\infty}^{\infty} da \, G((B_s)_{s \leq \tau_l^a})$$
By using invariance properties of Brownian motion for time and
space reversals, one obtains :
$$\int_0^{\infty} dt \, \mathbf{E} [ G((B_s)_{s \leq t})]= \int_0^{\infty}
dl \int_{- \infty}^{\infty} da \, \mathbf{E} [G((Z_s^{l,a})_{s \leq
  \tau_l + T_{a \rightarrow 0}})]$$
where $(Z_s^{l,a})_{s
\leq \tau_l + T_{a \rightarrow 0}}$ denotes a process such that $Z_s^{l,a} =
B_s$ 
for $s \leq
\tau_l$ and $(Z_{\tau_l + u}^{l,a})_{u \leq T_{a
    \rightarrow 0}}$ is the time-reversed process of a Brownian motion
starting from $a$, independent of $B$, and considered up to its first
hitting time of zero (denoted by $ T_{a
    \rightarrow 0}$).   \\ \\
Therefore, for all Borel sets $U$ of $\mathbf{R}_+^*$, if we define
$J_{c,U} (F)$ by : $$J_{c,U} (F) = \int_U dt \, \mathbf{E} \left[ F((L_t^y)_{y \in \mathbf{R}})
  \mathbf{1}_{|B_t| \geq c} \phi \left( \frac{1}{t} \int_{-c}^c L_t^y
  dy \right) \right]$$ we have, by taking $G_0$ and $G$ such that
$G((B_s)_{s \leq t}) = F((L_t^y)_{y \in \mathbf{R}})$ : 
$$ J_{c,U} (F) = \int_0^{\infty} dt \, \mathbf{E} \left[ F((L_t^y)_{y \in
    \mathbf{R}})  \mathbf{1}_{|B_t| \geq c} \phi \left( \frac{\int_{-c}^c L_t^y
  dy}{\int_{-\infty}^{\infty} L_t^y dy} \right)
    \mathbf{1}_{\int_{-\infty}^{\infty} L_t^y dy \in U} \right] $$
$$ = \int_0^{\infty} dl \int_{\mathbf{R} \backslash [-c,c]} da \,
 \mathbf{E} \left[ F((L^{y,l,a})_{y \in
    \mathbf{R}}) \phi \left(
 \frac{\int_{-c}^c L^{y,l,a}
  dy}{\int_{-\infty}^{\infty}  L^{y,l,a} dy} \right)
 \mathbf{1}_{\int_{-\infty}^{\infty}  L^{y,l,a} dy \in U} \right]$$
where  $(L^{y,l,a})_{y \in \mathbf{R}}$ is the continuous family of 
the total local times of $Z^{l,a}$. \\ \\
Hence, by Ray-Knight theorem applied to the independent processes \\ 
$(B_s=Z_s)_{s \leq \tau_l}$ and $(Z_{\tau_l + u})_{u \leq
 T_{a \rightarrow 0}}$, and classical additivity properties of squared
 Bessel processes : 
$$ J_{c,U} (F)=  \int_0^{\infty} dl \int_{\mathbf{R} \backslash
  [-c,c]} da \,  \mathbf{E} \left[ F((Y_{l,a}^y)_{y \in
    \mathbf{R}}) \phi \left(
 \frac{\mathcal{I}_{l,a}^c}{\mathcal{I}_{l,a}^{\infty}  } \right)
 \mathbf{1}_{\mathcal{I}_{l,a}^{\infty}  \in U} \right]   $$

$$ =  \int_0^{\infty} dl \int_{\mathbf{R} \backslash
  [-c,c]} da \,   \mathbf{E} \left[ F((Y_{l,a}^y)_{y \in
    \mathbf{R}}) \mathbf{E} \left[ \phi \left(
 \frac{\mathcal{I}_{l,a}^c}{\mathcal{I}_{l,a}^{\infty}  } \right)
 \mathbf{1}_{\mathcal{I}_{l,a}^{\infty}  \in U} \left| (Y_{l,a}^y)_{y \in
  [-c,c]} \right. \right] \right]$$
since $ F((Y_{l,a}^y)_{y \in \mathbf{R}})$ depends only on $ (Y_{l,a}^y)_{y \in
  [-c,c]}$. \\ \\
Now, if  $\theta$ is a given
  continuous function from $[-c,c]$ to $\mathbf{R}_+$, the integrals : $\int_c^{\infty} Y_{l,a}^y dy$ and
  $\int_{-\infty}^{-c} Y_{l,a}^y dy$ are independent conditionally on
  $(Y_{l,a}^y = \theta^y)_{y \in [-c,c]}$ and their
  conditional laws are respectively equal to the laws of
  $\int_0^{\infty} Y_{\theta^c, (a-c)_{+}}^y dy$ and 
 $\int_0^{\infty} Y_{\theta^{-c}, (-a-c)_{+}}^y dy$. \\ \\
Therefore, by additivity properties of BESQ processes, the
  conditional law of : $$\mathcal{I}_{l,a}^{\infty} -
  \mathcal{I}_{l,a}^{c} = \int_{- \infty}^{-c} Y_{l,a}^y dy +
  \int_{c}^{\infty} Y_{l,a}^y dy$$ given  $(Y_{l,a}^y = \theta^y)_{y
  \in [-c,c]}$,  is equal to the law of : $$
  \int_0^{\infty} Y_{\theta^c + \theta^{-c}, 0}^y dy + 
 \int_0^{\infty} Y_{0, (|a|-c)_+}^y dy$$ 
where $(Y_{\theta^c + \theta^{-c}, 0}^y)_{y \geq 0}$ and $(Y_{0,
  (|a|-c)_+}^y)_{y \geq 0}$ are supposed to be independent. \\ \\
By Ray-Knight theorem, $ \int_0^{\infty} Y_{\theta^c +
  \theta^{-c}, 0}^y dy$ has the same law as the time spent in
  $\mathbf{R}_+$ by $(B_s)_{s \leq \tau_{ \theta^c + \theta^{-c}}}$, therefore : $$ \int_0^{\infty} Y_{\theta^c +
  \theta^{-c}, 0}^y dy \overset{(d)}{=} \tau_{(\theta^c +
  \theta^{-c})/2} \overset{(d)}{=} T_{(\theta^c +
  \theta^{-c})/2}$$
Moreover : $$ \int_0^{\infty} Y_{0, (|a|-c)_+}^y dy \overset{(d)}{=}
  T_{(|a|-c)_+}$$
Hence, the conditional law of $\mathcal{I}_{l,a}^{\infty} -
  \mathcal{I}_{l,a}^{c}$, given  $(Y_{l,a}^y = \theta^y)_{y \in
  [-c,c]}$, is equal to the law of $T_{(|a|-c)_+ +(\theta^c +
  \theta^{-c})/2}$. 
Consequently : $$J_{c,U} (F) = \int_0^{\infty} dl \int_{\mathbf{R} \backslash
  [-c,c]} da \,  \mathbf{E} \left[ F((Y_{l,a}^y)_{y \in
    \mathbf{R}}) \psi_a(\mathcal{I}_{l,a}^c,\mathcal{Y}_{l,a}^c)
  \right]$$
where, for $|a| >c$ : 
$$\psi_a(\mathcal{I}, \theta) = \mathbf{E} \left[ \phi \left(
    \frac{\mathcal{I}}{\mathcal{I} + T_{|a|-c+\theta}} \right)
    \mathbf{1}_{\mathcal{I} + T_{|a|-c+\theta}  \in U} \right]$$
\noindent
Now, if, for all $u>0$, $p_u$ denotes the density of the law of $T_u$, one has : \\
$$\psi_a(\mathcal{I}, \theta) = \int_U \phi (\mathcal{I}/t) p_{|a|-c+
  \theta} (t-\mathcal{I}) dt$$
and : $$ J_{c,U} (F) =\int_U dt  \int_0^{\infty} dl \int_{\mathbf{R} \backslash
  [-c,c]} da \, \mathbf{E}  \left[  F((Y_{l,a}^y)_{y \in
    \mathbf{R}})  \phi \left( \frac{\mathcal{I}_{l,a}^c}{t} \right)
  p_{|a|-c+ \mathcal{Y}_{l,a}^c} (t-\mathcal{I}_{l,a}^c) \right]$$
By hypothesis, $  F((Y_{l,a}^y)_{y \in
    \mathbf{R}})$ depends only on $(Y_{l,a}^y)_{y \in
    [-c,c]}$. 
Moreover, for $a \geq c$, $(Y_{l,a}^y)_{y \in [-c,c]}$ has the same
  law as  $(Y_{l,+}^y)_{y \in [-c,c]}$, and for $a \leq -c$,
  $(Y_{l,a}^y)_{y \in [-c,c]}$ has the same law as
  $(Y_{l,-}^y)_{y \in [-c,c]}$. \\ \\
Hence, we have :
 $$ J_{c,U} (F) = \int_U dt  \int_0^{\infty} dl \, \mathbf{E}  \left[
  F((Y_{l,+}^y)_{y \in   \mathbf{R}}) \phi
  \left(\frac{\mathcal{I}_{l,+}^c}{t} \right) \int_c^{\infty} p_{a-c+
  \mathcal{Y}_{l,+}^c}(t-\mathcal{I}_{l,+}^c) da \right] $$
$$ +  \int_U dt  \int_0^{\infty} dl  \, \mathbf{E}  \left[
  F((Y_{l,-}^y)_{y \in   \mathbf{R}}) \phi
  \left(\frac{\mathcal{I}_{l,-}^c}{t} \right) \int_{-\infty}^{-c} p_{|a|-c+
  \mathcal{Y}_{l,-}^c}(t-\mathcal{I}_{l,-}^c) da \right]$$
Now, for $\theta \geq 0$, $u>0$ : 
$$  \int_{-\infty}^{-c} p_{|a|-c+ \theta} (u) da = \int_{c}^{\infty}
p_{a-c+ \theta} (u) da = \int_{\theta}^{\infty} p_b(u) db $$ $$= 
 \int_{\theta}^{\infty} \frac{b}{\sqrt{2 \pi u^3}} e^{-b^2/2u} db =
 \frac{1}{\sqrt{2 \pi u}} e^{-\theta^2/2u}$$
Therefore : $$  J_{c,U} (F) =  \int_U dt \, \frac{I_{c,t} (F)} {\sqrt{2
    \pi t}}$$ 
This equality is satisfied for every Borel set $U$. Hence, by definition
of $ J_{c,U} (F)$, the equality given in Proposition 2 occurs for
almost every $t > 0$. \\
In order to prove it for all $t> 0$, we begin to suppose that $F$ is
bounded and continuous. \\ 
In this case, for all $s$, $t >0$ :  
$$\left| \mathbf{E} \left[ F( (L_t^y)_{y \in \mathbf{R}}) \mathbf{1}_{|X_t|
    \geq c} \phi \left( \frac{1}{t} \int_{-c}^c L_t^y dy \right)
    \right] -\mathbf{E} \left[ F( (L_s^y)_{y \in \mathbf{R}}) \mathbf{1}_{|X_s|
    \geq c} \phi \left( \frac{1}{s} \int_{-c}^c L_s^y dy \right)
    \right] \right| $$ 
$$ \leq \mathbf{E}  \left[ \left| F( (L_t^y)_{y \in \mathbf{R}}) \phi \left(
    \frac{1}{t} \int_{-c}^c L_t^y dy \right) - F( (L_s^y)_{y \in \mathbf{R}})
 \phi \left( \frac{1}{s} \int_{-c}^c L_s^y dy \right)  \right| \right] $$ $$+
    ||F||_{\infty} \mathbf{P} (\exists u \in [s,t], |X_u|=c)$$
If $t$ is fixed, the first term of this sum tends to zero when $s$
    tends to $t$, by continuity of $F$, $\phi$ and dominated
    convergence. \\ 
The second term tends also to : $$ ||F||_{\infty} \mathbf{P} (
    |X_t|=c) = 0$$
Therefore, the function : $$t \rightarrow \mathbf{E} \left[ F( (L_t^y)_{y \in \mathbf{R}})
    \mathbf{1}_{|X_t|
   \geq c} \phi \left( \frac{1}{t} \int_{-c}^c L_t^y dy \right)
    \right]$$ is continuous. \\ \\
Now, let us prove that $I_{c,t}(F)$ is also continuous with respect to
    $t$. \\ \\
For all $t>0$ : $$ F((Y_{l,+}^y)_{y \in   \mathbf{R}})
    \frac{e^{-(\mathcal{Y}_{l,+}^c)^2/2(s -
    \mathcal{I}_{l,+}^c)}}{\sqrt{1 - \mathcal{I}_{l,+}^c/s}} \phi
    \left(\frac{\mathcal{I}_{l,+}^c}{s} \right) \underset{s
    \rightarrow t}{\rightarrow}   F((Y_{l,+}^y)_{y \in   \mathbf{R}})
    \frac{e^{-(\mathcal{Y}_{l,+}^c)^2/2(t -
    \mathcal{I}_{l,+}^c)}}{\sqrt{1 - \mathcal{I}_{l,+}^c/t}} \phi
    \left(\frac{\mathcal{I}_{l,+}^c}{t} \right) $$
by continuity of $\phi$ (if $\mathcal{I}_{l,+}^c < t$, it is
    clear, and if  $\mathcal{I}_{l,+}^c \geq t$, the two expressions
    are equal to zero for $s \leq 3t/2$). \\ \\
Moreover, for $s \leq 2t$ : 
$$ F((Y_{l,+}^y)_{y \in   \mathbf{R}})
    \frac{e^{-(\mathcal{Y}_{l,+}^c)^2/2(s -
    \mathcal{I}_{l,+}^c)}}{\sqrt{1 - \mathcal{I}_{l,+}^c/s}} \phi
    \left(\frac{\mathcal{I}_{l,+}^c}{s} \right) \leq \sqrt{3}
    ||F||_{\infty} e^{-(\mathcal{Y}_{l,+}^c)^2/4t} \leq \sqrt{3}
    ||F||_{\infty} e^{-(Y_{l,+}^c)^2/16t}$$
Recalling that the Lebesgue measure is invariant for the BESQ(2) process
    $(Y_{l,+}^y)_{y \geq 0}$, we have : 
$$\int_0^{\infty} dl \, \mathbf{E} \left[ e^{-(Y_{l,+}^c)^2/16t} \right] =
\int_0^{\infty} dl \, e^{-l^2/16t} < \infty$$
By dominated convergence, $t \rightarrow I_{c,t,+} (F)$ is
continuous. \\ \\
Similar computations imply the continuity of $t \rightarrow I_{c,t,-}
(F)$, and finally  $t \rightarrow I_{c,t}(F)$ is continuous. \\ \\
Consequently, for $F$ continuous and bounded, the equality given in
Proposition 2, which was proven for a.e. $t>0$, remains
true for every $t>0$.  \\ 
Now, by monotone class theorem (see
\cite{4}), it is not
difficult to extend this equality to every measurable and positive function,
which completes the proof of Proposition 2. \hfill$\Box$ \\ \\ 
This proposition has the following consequence : \\ \\
\textbf{Corollary 2 : } \textit{Let $F$ be a functional which satisfies the
condition of Proposition 2. The two following properties hold : \\ \\ 
1) For all $t>0$ : $$ \sqrt{2 \pi t} \, \mathbf{E} \left[ F((L_t^y)_{y
    \in \mathbf{R}})  \mathbf{1}_{|B_t| \geq c} \phi \left(
    \frac{1}{t} \int_{-c}^c L_t^y dy \right) \right] \leq \sqrt{3}
\, I(F) $$
2) When $t$ goes to infinity : $$ \sqrt{2 \pi t} \, \mathbf{E} \left[
  F((L_t^y)_{y   \in \mathbf{R}})  \mathbf{1}_{|B_t| \geq c} \phi \left(
    \frac{1}{t} \int_{-c}^c L_t^y dy \right) \right] \rightarrow
I(F)$$}
\textbf{Proof : } The first property is obvious, since $\phi(x) /
\sqrt{1-x} \leq \sqrt{3}$ for all $x \geq 0$. \\ 
In order to prove the second property, we distinguish two cases : \\
\\
1) If $I(F) < \infty$, we observe that : $$  F((Y_{l,+}^y)_{y \in
  \mathbf{R}})    \frac{e^{-(\mathcal{Y}_{l,+}^c)^2/2(t -
    \mathcal{I}_{l,+}^c)}}{\sqrt{1 - \mathcal{I}_{l,+}^c/t}} \phi
    \left(\frac{\mathcal{I}_{l,+}^c}{t} \right) $$
is smaller than $\sqrt{3}  F((Y_{l,+}^y)_{y \in
  \mathbf{R}})$ and tends to $ F((Y_{l,+}^y)_{y \in
  \mathbf{R}})$ when $t$ goes to infinity. \\ 
By dominated convergence, $I_{c,t,+}(F) \rightarrow I_+ (F)$. \\
Similarly,  $I_{c,t,-}(F) \rightarrow I_- (F)$ and finally :
 $$I_{c,t}(F) \rightarrow I (F)$$
2) If $I(F)=\infty$, we can suppose for example : $I_+(F) =
\infty$. \\
In this case : $$I_{c,t}(F) \geq I_{c,t,+}(F) \geq \int_0^{\infty} dl 
\, \mathbf{E} \left[  F((Y_{l,+}^y)_{y \in
  \mathbf{R}})   e^{-(\mathcal{Y}_{l,+}^c)^2/2(t -
    \mathcal{I}_{l,+}^c)} \phi
    \left(\frac{\mathcal{I}_{l,+}^c}{t} \right) \right]$$ 
which tends to $I_+(F)= \infty$ when $t \rightarrow
\infty$, by monotone convergence. \hfill$\Box$  \\ \\ 
Now, the next step in this article is the majorization of the
difference between the quantity $ \sqrt{2 \pi t} \, \mathbf{E}
[F((L_t^y)_{y \in \mathbf{R}})]$ and the expression given in
Proposition 2. 
\section{Majorization of the error term}
\noindent
For every positive and measurable functional $F$, we denote by
$\Delta_{c,t} (F)$ the error term we need to majorize : 
 $$\Delta_{c,t} (F) = \left| \sqrt{2 \pi t} \, \mathbf{E} \left[ F((L_t^y)_{y
    \in \mathbf{R}})  \mathbf{1}_{|B_t| \geq c} \phi \left(
    \frac{1}{t} \int_{-c}^c L_t^y dy \right) \right] -  \sqrt{2 \pi t} 
   \,  \mathbf{E} \left[ F((L_t^y)_{y
    \in \mathbf{R}}) \right] \right| $$
\noindent
It is easy to check that : 
$$ \Delta_{c,t}  (F) \leq  \Delta_{c,t}^{(1)}  (F) +  \Delta_{c,t}^{(2)}  (F)$$
where : $$ \Delta_{c,t}^{(1)}  (F) =  \sqrt{2 \pi t} 
\, \mathbf{E} \left[ F((L_t^y)_{y    \in \mathbf{R}})  \mathbf{1}_{|B_t|
     \leq c} \right]$$
and $$ \Delta_{c,t}^{(2)}  (F) =  \sqrt{2 \pi t} 
 \, \mathbf{E} \left[ F((L_t^y)_{y    \in \mathbf{R}})  \mathbf{1}_{ 
\int_{-c}^c L_t^y dy \geq t/3} \right]$$
The following proposition gives some precise majorizations of these
quantities, when $F$ satisfies the conditions of Definition 1. \\ \\
\textbf{Proposition 3 : }\textit{ Let $F$ be a functional from $\mathcal{C}
  (\mathbf{R}, \mathbf{R}_+ )$ to $ \mathbf{R}_+$ which satisfies the
  condition $C(c,n,h)$ for a positive, decreasing function $h$ and
  $c, n \geq 0$. \\
For all $t \geq 0$, one has the following majorizations : \\ \\
1) $\Delta_{c,t}^{(1)} (F) \leq A_n \frac{N_c(h)}{1 + (t/c^2)^{1/3}}$
  \\ \\
2)  $\Delta_{c,t}^{(2)} (F) \leq A_n \frac{ch(0)}{1 + (t/c^2)} \leq
  A_n \frac{N_c(h)}{1 + (t/c^2)}$ \\ \\
3)  $\Delta_{c,t}(F) \leq A_n \frac{N_c(h)}{1 + (t/c^2)^{1/3}}$
\\ \\
4) $I(F) \leq A_n N_c(h)$ \\ \\
where $A_n >0$ depends only on $n$.} \\ \\ 
In order to prove Proposition 3, we will need some inequalities  
 about the  processes
$(L_t^y)_{y \in [-c,c]}$ and $(Y_{l,+}^y)_{y \in [-c,c]}$. \\ 
More precisely, if we put : $\Sigma_t^c = \underset{y \in
  [-c,c]}{\sup} L_t^y$,  $\sigma_t^c = \underset{y \in
  [-c,c]}{\inf} L_t^y$, $\Theta_{l,+}^c = \underset{y \in
  [-c,c]}{\sup} Y_{l,+}^y$,  $\theta_{l,+}^c = \underset{y \in
  [-c,c]}{\inf} Y_{l,+}^y$, $\Theta_{l,-}^c = \underset{y \in
  [-c,c]}{\sup} Y_{l,-}^y$,  $\theta_{l,-}^c = \underset{y \in
  [-c,c]}{\inf} Y_{l,-}^y$, the following statement hold : \\ \\
\textbf{Lemma 3 : } \textit{For all $c$, $t > 0$ : \\  \\
1) If $a \geq 0$ :  $$\mathbf{P} \left( \frac {\Sigma_t^c + c}
  {\sigma_t^c + c} \geq a \right) \leq A e^{-\lambda a}$$
2) If $a \geq 4$ :  $$\mathbf{P} \left( \frac {\Theta_{l,+}^c + c}
  {\theta_{l,+}^c + c} \geq a \right) \leq A e^{-\lambda \left(a+
    \frac{l}{c} \right)}$$
3)  If $a \geq 4$ :  $$\mathbf{P} \left( \frac {\Theta_{l,-}^c + c}
  {\theta_{l,-}^c + c} \geq a \right) \leq A e^{-\lambda \left(a+
    \frac{l}{c} \right)}$$   
where $A>0$, $0<\lambda<1$ are universal constants.} \\ \\
\textbf{Proof of Lemma 3 : } 
 1) Let us suppose $a \geq 8$, $c > 0$. \\
 In that case : $$ \mathbf{P} \left( \frac {\Sigma_t^c + c}
  {\sigma_t^c + c} \geq a, L_t^0 \geq \frac{ac}{4} \right) \leq 
 \mathbf{P} \left( \frac {\Sigma_t^c + c} {\sigma_t^c + c} \geq
 8 , L_t^0 \geq \frac{ac}{4}  \right)$$
$$ \leq  \underset{k \in \mathbf{N}}{\sum} \mathbf{P} \left(
  \frac{\Sigma_t^c}{\sigma_t^c} \geq 8, L_t^0 \in [2^{k-2} ac, 2^{k-1}
  ac] \right) $$ $$  \leq  \underset{k \in \mathbf{N}}{\sum}
  \mathbf{P} ( \Sigma_t^c \geq 2^{k} ac, L_t^0 \in [2^{k-2} ac, 2^{k-1}
  ac] ) $$ $$ +   \underset{k \in \mathbf{N}}{\sum}  \mathbf{P} ( \sigma_t^c \leq 2^{k-3} ac,  \Sigma_t^c \leq
  2^{k} ac,  L_t^0 \in [2^{k-2} ac, 2^{k-1}
  ac] ) $$ $$\leq  \underset{k \in \mathbf{N}}{\sum} \left[
  \mathbf{P} ( \Sigma_{\tau_{2^{k-1} ac}}^c \geq 2^k ac) + 
 \mathbf{P} ( \sigma_{\tau_{2^{k-2} ac}}^c \leq 2^{k-3} ac, \,
 \Sigma_{\tau_{2^{k-2} ac}}^c \leq 2^{k} ac ) \right] $$ $$= 
 \underset{k \in \mathbf{N}}{\sum} \left[ \alpha_c (2^{k-1} a c) +
  \beta_c (2^{k-2} ac) \right] $$
where for $l \geq 0$, $\alpha_c (l) = \mathbf{P} (\Sigma_{\tau_l}^c \geq
  2 l)$ and  $\beta_c (l) = \mathbf{P} (\sigma_{\tau_l}^c \leq
  l/2, \Sigma_{\tau_l}^c \leq
  4l )$. \\ \\
Now, by Ray-Knight theorem, $\alpha_c (l) \leq 2 \mathbf{P} \left(
  \underset{y \in [0,c]}{\sup} Y_{l,0}^y \geq 2l \right)$, and by 
Dubins-Schwarz theorem, $Y_{l,0}^y = l + \beta_{\int_0^y 4 Y_{l,0}^z
  dz} $, where $\beta$ is a Brownian motion. \\ \\
Hence, if $S = \inf \{ y \geq 0, Y_{l,0}^y \geq 2l \}$, one has :  
$ \underset{u \leq \int_0^S 4 Y_{l,0}^z dz}{\sup} \beta_u = l $, 
and if we suppose $ \underset{y \in [0,c]}{\sup} Y_{l,0}^y \geq
2l$, we have $S \leq c$,  $\int_0^S 4 Y_{l,0}^z dz \leq \int_0^S 8l dz
\leq 8lc$, and finally : $\underset{u \leq 8lc}{\sup} \beta_u \geq
l$. \\ \\
Consequently : $$ \alpha_c (l) \leq 2 \mathbf{P} \left( \underset{u
  \leq 8lc}{\sup} \beta_u \geq l \right) = 2 \mathbf{P} (|\beta_{8lc}|
  \geq l) \leq 4 \mathbf{P} (\beta_{8lc} \geq l) \leq 4 e^{-l/16c}$$
By the same kind of argument, one obtains : $$\beta_c (l) \leq 4
e^{-l/128c}$$ and finally : $$ \mathbf{P} \left( \frac {\Sigma_t^c + c}
  {\sigma_t^c + c} \geq a, L_t^0 \geq \frac{ac}{4} \right) \leq 4
\underset{k \in \mathbf{N}}{\sum} \left( e^{-2^{k-1} a/16} +
  e^{-2^{k-2} a /128} \right) $$ $$\leq 8 \underset{k \in
  \mathbf{N}}{\sum} e^{-2^k a/512} \leq 8 \underset{k \in
  \mathbf{N^*}}{\sum} e^{-k a/512} \leq 8 e^{-a/512} \left( \underset{k \in
  \mathbf{N}}{\sum} e^{-k/64} \right) \leq 520 e^{-a/512}$$
On the other hand : $$ \mathbf{P} \left( \frac {\Sigma_t^c + c}
  {\sigma_t^c + c} \geq a, L_t^0 \leq \frac{ac}{4} \right) \leq
  \mathbf{P}  \left(\Sigma_t^c + c \geq ac,  L_t^0 \leq \frac{ac}{4}
  \right) \leq  \mathbf{P}  \left(\Sigma_{\tau_{ac/4}}^c \geq (a-1)c
  \right) $$ $$  \leq  \mathbf{P}  \left(\Sigma_{\tau_{ac/4}}^c \geq
  \frac{7ac}{8}  \right) \leq \alpha_c \left(\frac{ac}{4}\right) \leq
  4 e^{-a/64}$$ 
Consequently : $$ \mathbf{P} \left( \frac {\Sigma_t^c + c}
  {\sigma_t^c + c} \geq a \right) \leq 524 e^{-a/512}$$ for all $a
  \geq 8$. \\ 
This inequality remains obviously true for $a \leq 8$ or $c=0$, so the
  first part of Lemma 3 is proven. \\ \\ \\
2) Let $a$ be greater than 4. If $l \geq ac/4$ :  
$$ \mathbf{P} \left( \frac{\Theta_{l,+}^c + c}{\theta_{l,+}^c + c}
  \geq 4 \right) \leq  \mathbf{P} \left(\Theta_{l,+}^c \geq 2l \right)
  +  \mathbf{P} \left(\Theta_{l,+}^c \leq 2l,\theta_{l,+}^c \leq l/2
  \right) \leq 2 \tilde{\alpha}_c(l) + \tilde{\beta}_c(l)$$
where  $$\tilde{\alpha}_c(l) =  \mathbf{P} \left(\underset{y \in
  [0,c]}{\sup} Y_{l,+}^y  \geq 2l \right)$$ and  $$\tilde{\beta}_c(l) =
  \mathbf{P} \left(\underset{y \in [-c,c]}{\sup} Y_{l,0}^y  \leq 2l, 
\underset{y \in [-c,c]}{\inf} Y_{l,0}^y  \leq l/2 \right)$$. \\ \\
Now, $(Y_{l,+}^y)_{y \geq 0}$ is a BESQ(2) process, hence, if
  $(\beta_y= (\beta_y^{(1)}, \beta_y^{(2)}))_{y \geq 0}$ is a standard
  two-dimensional Brownian motion : 
$$ \tilde{\alpha}_c(l) = \mathbf{P} \left( \underset{y \in
  [0,c]}{\sup}  Y_{l,+}^y  \geq 2l \right) =  \mathbf{P}  \left(
  \underset{y \leq c}{\sup} \, ||\beta_y + (\sqrt{l},0)|| \geq \sqrt{2l}
  \right) $$ $$ \leq  \mathbf{P} \left( \underset{y \leq c}{\sup}
 \, ||\beta_y|| \geq \sqrt{l} (\sqrt{2} - 1) \right) \leq 2   \mathbf{P}
  \left( \underset{y \leq c}{\sup} \, | \beta_y^{(1)}| \geq \sqrt{l}
  \left( \frac{\sqrt{2}-1}{2} \right)\right) $$ $$ \leq 8  \mathbf{P}
  \left(\beta_c^{(1)} \geq \sqrt{l}
  \left( \frac{\sqrt{2}-1}{2} \right) \right) \leq 8e^{-l/50c}$$
Moreover : $$ \tilde{\beta}_c(l) \leq \mathbf{P}  \left( \underset{y
  \in [-c,c]}{\sup} Y_{l,0}^y  \leq 4l,  \underset{y \in [-c,c]}
{\inf} Y_{l,0}^y  \leq l/2 \right) = \beta_c(l) \leq 4 e^{-l/128 c}$$ 
Therefore, if $l \geq ac/4$ : $$\mathbf{P} \left( \frac{\Theta_{l,+}^c
  + c}{\theta_{l,+}^c + c}  \geq a \right) \leq 20 e^{-l/128c}$$ 
Now, let us suppose  $l \leq ac/4$. In this case : 
$$\mathbf{P} \left( \frac{\Theta_{l,+}^c
  + c}{\theta_{l,+}^c + c}  \geq a \right) \leq \mathbf{P}  \left(
\Theta_{ac/4,+}^c \geq 3ac/4 \right) \leq 2 \tilde{\alpha}_c(ac/4)
  \leq 16 e^{-a/200}$$
Hence, for every $l \geq 0$, $a \geq 4$ : $$ \mathbf{P} \left(
  \frac{\Theta_{l,+}^c + c}{\theta_{l,+}^c + c}  \geq a \right) \leq
  20 e^{-(a + (l/c))/1024}$$ which proves the second inequality of the
  lemma. \\ The proof of the third inequality is exactly similar.
  \hfill$\Box$  \\
  \\
Now, we are able to prove the main result of the section, which was
  presented in Proposition 3. \\ \\
\textbf{Proof of Proposition 3 : } 1)  For $c=0$,  $\Delta_{c,t}^{(1)} (F) =0$, so we can
  suppose $c >0$. \\ \\
The functional $F$ satisfies the condition $C(c,n,h)$; hence, for all
 $a \geq 1$ : $$ \frac{\Delta_{c,t}^{(1)} (F)}{\sqrt{2 \pi t}}
  = \mathbf{E} \left[ F((L_t^y)_{y \in \mathbf{R}}) \mathbf{1}_{|B_t|
  \leq c} \right] $$ $$ \leq \mathbf{E}  \left[ \left(
  \frac{\Sigma_t^c + c}{\sigma_t^c + c} \right) ^n h(\sigma_t^c)
  \mathbf{1}_{|B_t|  \leq c} \right] $$ $$\leq  \mathbf{E}  \left[
 \left(\frac{\Sigma_t^c + c}{\sigma_t^c + c} \right) ^n  h(0)
  \mathbf{1}_{\frac{\Sigma_t^c + c}{\sigma_t^c + c} \geq a} \right] +
  a^n \mathbf{E}  \left[ h(\sigma_t^c) \mathbf{1}_{|B_t|  \leq c}
  \mathbf{1}_{\frac{\Sigma_t^c + c}{\sigma_t^c + c} \leq a} \right]$$
Now, if $\frac{\Sigma_t^c + c}{\sigma_t^c + c} \leq a$, 
 $\frac{L_t^0 + c}{\sigma_t^c + c} \leq a$ and $\sigma_t^c \geq
  \left( \frac{L_t^0}{a} - c \right)_+ $. \\ \\
Therefore : $$  \frac{\Delta_{c,t}^{(1)} (F)}{\sqrt{2 \pi t}}  \leq h(0) 
 \mathbf{E}  \left[ \left( \frac{\Sigma_t^c + c}{\sigma_t^c + c}
  \right) ^n  \mathbf{1}_{\frac{\Sigma_t^c + c}{\sigma_t^c + c} \geq
  a} \right] + a^n  \mathbf{E} \left[ h \left(  \left(
  \frac{L_t^0}{a} - c \right)_+ \right)  \mathbf{1}_{|B_t|  \leq c}
  \right]$$
By Lemma 3 : $$   \mathbf{E}  \left[ \left( \frac{\Sigma_t^c +
  c}{\sigma_t^c + c}   \right) ^n  \mathbf{1}_{\frac{\Sigma_t^c +
  c}{\sigma_t^c + c} \geq  a} \right] = a^n \mathbf{P} \left(
  \frac{\Sigma_t^c +  c}{\sigma_t^c + c} \geq  a \right) +
  \int_{a}^{\infty} n b^{n-1}  \mathbf{P} \left(
  \frac{\Sigma_t^c +  c}{\sigma_t^c + c} \geq  b \right) db $$
$$ \leq A \left( a^n e^{-\lambda a} +  \int_{a}^{\infty} n b^{n-1}
e^{-\lambda b} db \right) = A a^n e^{- \lambda a} \left( 1 + n \int_{0}^{\infty}
  \frac{(a+b)^{n-1}}{a^n} e^{-\lambda b} db \right)$$ $$ \leq  A a^n
e^{- \lambda a} \left( 1 + n \int_{0}^{\infty}
 (1+b)^{n} e^{-\lambda b} db \right) \leq A \left( \frac{6}{\lambda}
\right) ^{n+1} (n+1)! \, a^n  e^{-\lambda a}$$
On the other hand, by using the probability density of $(L_t^0,
|B_t|)$ (given for example in \cite{3}, Lemma 2.4) : $$  \mathbf{E}  \left[ h \left(  \left(
  \frac{L_t^0}{a} - c \right)_+ \right)  \mathbf{1}_{|B_t|  \leq c}
  \right] $$ $$ = \sqrt{\frac{2}{\pi t^3}} \int_0^{\infty} dl \int_0^{c} dx
\,  h \left(  \left(  \frac{l}{a} - c \right)_+ \right) (l+x)
  e^{-(l+x)^2/2t} $$ $$= \sqrt{\frac{2}{\pi t^3}} h(0)  \int_0^{ac} dl
 \int_0^{c} dx \, (l+x)  e^{-(l+x)^2/2t} $$ $$ +  \sqrt{\frac{2}{\pi t^3}} 
 \int_{ac}^{\infty} dl  \int_0^{c} dx \, h \left(  \frac{l}{a} - c
 \right) (l+x)  e^{-(l+x)^2/2t}$$ $$ =  \sqrt{\frac{2}{\pi}}
 \frac{c^2}{t} h(0) \int_0^{a} dl \int_0^{1} dx \,
 \frac{c(l+x)}{\sqrt{t}} e^{-c^2 (l+x)^2/2t}$$ $$  +  \sqrt{\frac{2}{\pi}}
 \frac{ac^2}{t} \int_0^{\infty} dl \int_0^{1} dx \, h(cl) 
 \frac{c(al+a+x)}{\sqrt{t}}  e^{-c^2 (al+a + x)^2/2t} $$
For all $\theta \geq 0$, $\theta e^{-\theta^2/2} \leq e^{-1/2} \leq
1$. Hence : $$  \mathbf{E}  \left[ h \left(  \left(
  \frac{L_t^0}{a} - c \right)_+ \right)  \mathbf{1}_{|B_t|  \leq c}
  \right] \leq  \sqrt{\frac{2}{\pi}} \frac{a c^2}{t} \left( h(0) +
    \int_0^{\infty} h(cl) dl \right) =  \sqrt{\frac{2}{\pi}} \frac{a
    c}{t} N_c(h)$$
Moreover, for $0<t \leq c^2$ : 
$$ \mathbf{E}  \left[ h \left(  \left(
  \frac{L_t^0}{a} - c \right)_+ \right)  \mathbf{1}_{|B_t|  \leq c}
  \right] \leq h(0) \leq \frac{N_c(h)}{c} \leq \frac{a
  N_c(h)}{\sqrt{t}}$$
The majorizations given above imply : $$ 
\Delta_{c,t}^{(1)} (F) \leq  A \left( \frac{6}{\lambda}
\right) ^{n+1} (n+1)! \, a^n  e^{-\lambda a} \sqrt{2 \pi t} \, h(0) + 
 \sqrt{2 \pi } \, a^{n+1} \left( \frac{c}{\sqrt{t}} \wedge 1 \right) N_c(h)$$
Now, let us choose $a$ as a function of $t$. \\ \\
For $t \leq c^2$, we take $a=1$ and obtain : 
$$ \Delta_{c,t}^{(1)} (F)  \leq A \left( \frac{6}{\lambda}
\right) ^{n+1} (n+1)! \,  e^{-\lambda} \sqrt{2 \pi} \, c h(0) +  \sqrt{2
  \pi} \,  N_c(h)$$ $$ \leq \sqrt{2 \pi} \left(  1 +  A \left( \frac{6}{\lambda}
\right) ^{n+1} (n+1)! \,  e^{-\lambda} \right) N_c(h) $$
For $t \geq c^2$, we take $a = (t/c^2)^{1/6(n+1)}$ : 
$$ \Delta_{c,t}^{(1)} (F)  \leq  A \left( \frac{6}{\lambda}
\right) ^{n+1} (n+1)! \left( \frac{t}{c^2} \right) ^{1/6} e^{- \lambda
  \left( \frac{t}{c^2} \right) ^{1/6(n+1)}} \sqrt{2 \pi t} \,  h(0) +
\sqrt{2 \pi}  \left( \frac{t}{c^2} \right) ^{1/6} \frac{c}{\sqrt{t}}
N_c(h)$$
$$\leq  \sqrt{2 \pi}  \left( 1+ A \left( \frac{6}{\lambda}
\right) ^{n+1} (n+1)! \right) N_c(h)  \left( \frac{t}{c^2} \right) ^{-1/3} 
  \left( 1 +  \frac{t}{c^2} e^{- \lambda
  \left( \frac{t}{c^2} \right) ^{1/6(n+1)}} \right)$$ 
$$\leq  \sqrt{2 \pi}  \left( 1+ A \left( \frac{6}{\lambda}
\right) ^{n+1} (n+1)! \right)   \left( 1 + \underset{u \geq 1}{\sup} \, u
e^{- \lambda u^{1/6(n+1)}}
 \right) \left( \frac{t}{c^2} \right)^{-1/3}  N_c(h) $$
where $ \underset{u \geq 1}{\sup} \, u e^{- \lambda u
    ^{1/6(n+1)}}$ is finite and depends only on $n$ (we recall the
  $\lambda$ is a universal constant). \\ \\
In the two cases, the first inequality of Proposition 3 is
satisfied. \\ \\
2) For $c=0$, $\Delta_{c,t}^{(2)} (F)= 0$, so we can again suppose
$c>0$. \\ 
For $a \geq 1$ : 
$$ \frac{\Delta_{c,t}^{(2)} (F)}{\sqrt{2 \pi t}} = \mathbf{E} \left[
  F((L_t^y)_{y \in \mathbf{R}} ) \mathbf{1}_{\int_{-c}^c L_t^y dy \geq
  t/3} \right] \leq \mathbf{E} \left[ \left( \frac{\Sigma_t^c +
  c}{\sigma_t^c+ c} \right) ^n h(\sigma_t^c)  \mathbf{1}_{\Sigma_t^c
  \geq t/6c} \right]  $$ $$ \leq h(0) \left(  \mathbf{E}  \left[ \left(
 \frac{\Sigma_t^c + c}{\sigma_t^c+ c} \right) ^n \mathbf{1}_{
  \frac{\Sigma_t^c + c}{\sigma_t^c+ c} \geq a} \right] + a^n
  \mathbf{P} \left( L_t^0 \geq \frac{t}{6ac} - c \right) \right) $$
 $$ \leq A \left( \frac{6}{\lambda}
\right) ^{n+1} (n+1)! \, a^n e^{- \lambda a} h(0) + 2 a^n h(0) e^{-
  \frac{1}{2t}\left( \frac{t}{6ac} -c \right)_+^2} $$
If $t \leq 12 c^2$, we take $a=1$ : 
$$  \Delta_{c,t}^{(2)} (F)\leq c h(0) \sqrt{24 \pi}
\left( 2 + 
 A \left( \frac{6}{\lambda} \right) ^{n+1} (n+1)! \, e^{-\lambda}
 \right)$$ 
 If $t \geq 12 c^2$, we take $a=\left( \frac{t}{12c^2} \right)^{1/3}$
 : 
$$ \Delta_{c,t}^{(2)} (F) \leq \sqrt{2 \pi t} \, h(0) \left[ A \left(
    \frac{6}{\lambda} \right) ^{n+1} (n+1)! \left(\frac{t}{12 c^2}
    \right)^{n/3} e^{- \lambda \left( \frac{t}{12 c^2} \right) ^{1/3}}
... \right.$$ $$ \left. ... + 2  \left(\frac{t}{12 c^2} \right)^{n/3} e^{- \frac{c^2}{2t}
    \left( 2 (t/12c^2) ^{2/3} - 1 \right) ^2 } \right]$$ 
$$ \leq \left(\frac{c^2}{t} \right) c h(0) \sqrt{2 \pi} \, 12^{3/2} \left( 2+ 
A \left(  \frac{6}{\lambda} \right) ^{n+1} (n+1)! \right)  \left(
\frac{t}{12c^2} \right) ^{\frac{n}{3} + \frac{3}{2}} \left(e^{-\lambda
 \left( \frac{t}{12c^2} \right) ^{1/3}} + e^{-\frac{1}{24}
  \left( \frac{t}{12c^2} \right) ^{1/3}} \right)$$
The second inequality of Proposition 3 holds, since $\underset{u \geq
  1}{\sup} \, u^{\frac{n}{3} + \frac{3}{2}} \left( e^{- \lambda u^{1/3}}
  +  e^{- \frac{1}{24}\lambda u^{1/3}} \right)$ is finite and depends only on
  $n$. \\ \\
3) This inequality is an immediate consequence of 1) and 2). \\ \\
4) For every $l \geq 0$ : $$\mathbf{E} [ F((Y_{l,+}^y)_{y \in
  \mathbf{R}})] \leq \mathbf{E} \left[ \left(
  \frac{\Theta_{l,+}^c+c}{\theta_{l,+}^c + c} \right)^n
h(\theta_{l,+}^c) \right] $$ $$ \leq h(0) \mathbf{E}  \left[ \left(
  \frac{\Theta_{l,+}^c+c}{\theta_{l,+}^c + c} \right)^n   \mathbf{1}_{
 \frac{\Theta_{l,+}^c+c}{\theta_{l,+}^c + c} \geq 4} \right] + 4^n h
\left( \left( \frac{l}{4} - c \right)_+ \right) $$ 
Now, by Lemma 3 : 
$$ \mathbf{E}  \left[ \left(
  \frac{\Theta_{l,+}^c+c}{\theta_{l,+}^c + c} \right)^n   \mathbf{1}_{
 \frac{\Theta_{l,+}^c+c}{\theta_{l,+}^c + c} \geq 4} \right] $$ $$ = 
4^n  \mathbf{P} \left(  \frac{\Theta_{l,+}^c+c}{\theta_{l,+}^c + c}
  \geq 4 \right)  + \int_{4}^{\infty} n b^{n-1}  \mathbf{P} \left(
  \frac{\Theta_{l,+}^c+c}{\theta_{l,+}^c + c} \geq b \right) db$$
$$ \leq A e^{-\lambda l/c} \left( 4^n e^{-4 \lambda} + \int_4^{\infty} n
b^{n-1} e^{-\lambda b} db \right)  $$ $$ \leq  A e^{-\lambda l/c}
\left(\frac{6}{\lambda} \right) ^{n+1} (n+1)! \, 4^n e^{-4 \lambda}$$
Hence : $$  \mathbf{E}  [ F((Y_{l,+}^y)_{y \in
  \mathbf{R}})] \leq A h(0)  e^{-\lambda l/c} \left(\frac{6}{\lambda}
\right) ^{n+1} (n+1)! \, 4^n e^{-4 \lambda} + 4^n  h
\left( \left( \frac{l}{4} - c \right)_+ \right)$$ and, by integrating
with respect to $l$ : 
$$I_+(F) \leq \frac{A}{\lambda}  \left(\frac{6}{\lambda} \right)
^{n+1} (n+1)! 4^n e^{-4 \lambda} c h(0) + 4^{n+1} c h(0) + 4^{n+1}
\int_0^{\infty} h(l) dl $$ $$ \leq 4^{n+1}  \left(1 +  \frac{A}{\lambda}
 \left(\frac{6}{\lambda} \right)
^{n+1} (n+1)! \right) N_c(h)  $$
By symmetry, the same inequality holds for $I_-(F)$, and : 
$$I(F) \leq  2^{2n+3}  \left(1 +  \frac{A}{\lambda}
 \left(\frac{6}{\lambda} \right)
^{n+1} (n+1)! \right) N_c(h)$$
which completes the proof of Proposition 3.  \hfill$\Box$
\section{An estimation of the quantity : $\mathbf{E} [F((L_t^y)_{y \in
 \mathbf{R}})]$}
\noindent
In this section, we majorize $\mathbf{E} [F((L_t^y)_{y \in
 \mathbf{R}})]$ by an equivalent of this quantity when $t$ goes to
 infinity. The following statement holds : \\ \\
\textbf{Proposition 4.1 : } \textit{ Let $F$ be a functional from
 $\mathcal{C}(\mathbf{R}, \mathbf{R}_+)$ to $\mathbf{R}_+$, which
 satisfies the condition $C(c, n,h)$, for a positive, decreasing
 function $h$, and $c$, $n \geq 0$. \\ 
The following properties hold : \\ \\
1) For all $t > 0$ : $$\sqrt{2 \pi t} \, \mathbf{E} [F ((L_t^y)_{y \in
 \mathbf{R}}) ] \leq K_n N_c(h)$$ 
where $K_n > 0$ depends only on $n$. \\ \\
2) If $N_c(h) < \infty$ : 
$$\sqrt{2 \pi t} \mathbf{E} [F ((L_t^y)_{y \in
 \mathbf{R}}) ] \underset{t \rightarrow \infty}{\rightarrow} I(F)$$} 
\noindent 
\textbf{Proof : } We suppose $N_c(h) < \infty$. \\ 
 Proposition 3 implies the following : $$ 
\Delta_{c,t} (F) \leq A_n N_c(h) $$
$$\Delta_{c,t} (F)  \underset{t
  \rightarrow \infty}{\rightarrow} 0$$
\noindent
Moreover, by Corollary 2 : 
$$ \sqrt{2 \pi t} \, \mathbf{E} \left[ F ((L_t^y)_{y \in
 \mathbf{R}})\mathbf{1}_{|B_t| \geq c} \phi \left( \frac{1}{t}
 \int_{-c}^c L_t^y dy \right) \right]   \underset{t
  \rightarrow \infty}{\rightarrow} I(F) $$
$$ \sqrt{2 \pi t} \, \mathbf{E} \left[ F ((L_t^y)_{y \in
 \mathbf{R}})\mathbf{1}_{|B_t| \geq c} \phi \left( \frac{1}{t}
 \int_{-c}^c L_t^y dy \right) \right]   \leq \sqrt{3} \, I(F) \leq
 \sqrt{3} \,  A_n N_c(h) $$
for all $t > 0$. \\ \\
Now, by definition, one has : 
$$\left| \sqrt{2 \pi t} \, \mathbf{E} [F ((L_t^y)_{y \in
 \mathbf{R}}) ] -  \sqrt{2 \pi t} \, \mathbf{E} \left[ F ((L_t^y)_{y \in
 \mathbf{R}})\mathbf{1}_{|B_t| \geq c} \phi \left( \frac{1}{t}
 \int_{-c}^c L_t^y dy \right) \right] \right| = \Delta_{c,t} (F)$$ 
Therefore : $$ \sqrt{2 \pi t} \, \mathbf{E} [F ((L_t^y)_{y \in
 \mathbf{R}}) ] \underset{t
  \rightarrow \infty}{\rightarrow} I(F)$$
$$ \sqrt{2 \pi t} \, \mathbf{E} [F ((L_t^y)_{y \in
 \mathbf{R}}) ] \leq (1 + \sqrt{3}) A_n N_c(h)$$
which proves Proposition 4.1.  \hfill$\Box$   \\ \\
The following result is an extension of Proposition 4.1 to a larger
 class of functionals $F$ : \\ \\
\textbf{Proposition 4.2 : } \textit{ Let $F : \mathcal{C}( \mathbf{R},
 \mathbf{R}_+) \rightarrow \mathbf{R}_+$ be a positive and measurable
 functional. The following properties hold for all $n \geq 0$ : \\ \\
1) For all $t > 0$ : 
$$ \sqrt{2 \pi t} \, \mathbf{E} [F((L_t^y)_{y \in \mathbf{R}})] \leq K_n
N^{(n)} (F) $$ 
\noindent
2) If $N^{(n)} (F) < \infty$ : $$\sqrt{2 \pi t} \, \mathbf{E}
   [F((L_t^y)_{y \in \mathbf{R}})] \underset{t \rightarrow
   \infty}{\rightarrow} I(F)$$} 
\noindent
\textbf{Proof : } We suppose $N^{(n)} (F) <\infty$. \\ \\
1) Let us take $M$ such that $N^{(n)} (F) <M$. \\ 
By definition, $F$ satisfies the condition $D(n, M)$, so there exists
$(c_k)_{k \geq 1}$,$(h_k)_{k \geq 1}$, $(F_k)_{k \geq 0}$ as in
Definition 2. \\ \\
One has : $F = \underset{k \geq 1}{\sum} (F_k - F_{k-1})$, hence : 
$$ \sqrt{2 \pi t} \, \mathbf{E}
   [F((L_t^y)_{y \in \mathbf{R}})] \leq   \underset{k \geq 1}{\sum} 
\sqrt{2 \pi t} \, \mathbf{E}
   [|F_k - F_{k-1}|((L_t^y)_{y \in \mathbf{R}})]$$ $$ \leq K_n 
 \underset{k \geq 1}{\sum} N_{c_k}(h_k) \leq K_n M$$
By taking $M \rightarrow N^{(n)} (F)$, one obtains the first part of
   Proposition 4.2. \\ \\
2) In order to prove the convergence, let us consider the equality :
   $$  \sqrt{2 \pi t} \, \mathbf{E}
   [F((L_t^y)_{y \in \mathbf{R}})] =  \underset{k \geq 1}{\sum}
\sqrt{2 \pi t} \,\mathbf{E}
   [(F_k - F_{k-1})_+ ((L_t^y)_{y \in \mathbf{R}})] $$ $$ -
   \underset{k \geq 1}
{\sum}
\sqrt{2 \pi t} \, \mathbf{E}
   [(F_k - F_{k-1})_- ((L_t^y)_{y \in \mathbf{R}})]$$
where the two sums are convergent. \\ \\
By Proposition 4.1, the two terms indexed by $k$ tend to $I((F_k -
   F_{k-1})_+)$ and $I((F_k - F_{k-1})_-)$ when $t$ goes to infinity,
   and they are bounded by $K_n N_{c_k} (h_k)$. \\ \\
Hence, by dominated convergence : 
$$ \sqrt{2 \pi t} \, \mathbf{E}
   [F((L_t^y)_{y \in \mathbf{R}})] \underset{t \rightarrow
   \infty}{\rightarrow} \underset{k \geq 1}{\sum} I((F_k -
   F_{k-1})_+) - \underset{k \geq 1}{\sum} I((F_k - F_{k-1})_-)$$
\noindent 
Now, by definition of $I$ : 
$$ \underset{k \geq 1}{\sum} I((F_k -
   F_{k-1})_+) = I \left(  \underset{k \geq 1}{\sum} (F_k - F_{k-1})_+
   \right)$$
$$ \underset{k \geq 1}{\sum} I((F_k -
   F_{k-1})_-) = I \left(  \underset{k \geq 1}{\sum} (F_k - F_{k-1})_-
   \right)$$
Therefore, if $G =  \underset{k \geq 1}{\sum} (F_k - F_{k-1})_+$, and 
$H =  \underset{k \geq 1}{\sum} (F_k - F_{k-1})_-$, one has : $$
 \underset{k \geq 1}{\sum} I((F_k -
   F_{k-1})_+) - \underset{k \geq 1}{\sum} I((F_k - F_{k-1})_-) = I(G)
   - I(H)$$ 
where : 
$$ I(G) - I(H) = I(G-H) = I(F)$$ 
since $I(G) + I(H) \leq 2 K_n  \underset{k \geq 1}{\sum} N_{c_k} (h_k)
< \infty$.  \hfill$\Box$  \\ \\
Proposition 4.2 is proven, and we now have all we need for the proof of the
main Theorem, which is given in Section 5.  
\section{Proof of the main Theorem}
\noindent
Our proof of the Theorem starts with a general lemma (which does not
involve Wiener measure) : \\ \\
\textbf{Lemma 5 : } \textit{ If $F : \mathcal{C} (\mathbf{R}, \mathbf{R}_+)
\rightarrow \mathbf{R}_+$ is a measurable functional, $l_0
\in \mathcal{C} (\mathbf{R}, \mathbf{R}_+)$, $x \in \mathbf{R}$, and
$n \geq 0$ :
 $$ N^{(n)} \left(  F^{(l_0^y)_{y \in \mathbf{R}}, x} \right) \leq 2^n \left( 1 +
 \left( \underset{z \in \mathbf{R}}{\sup} \, l_0^z \right)^n \right)
 \left( 1 + |x| \right)^{n+1} N^{(n)} (F) $$}
\noindent
\textbf{Proof of Lemma 5 : } Let $M$ be greater than $N^{(n)} (F)$. \\ 
There exists a sequence $(c_k)_{k \geq 1}$ in $[1, \infty[$, a
sequence $(h_k)_{k \geq 1}$ of decreasing functions from
$\mathbf{R}_+$ to $\mathbf{R}_+$, and a sequence $(F_k)_{k \geq 0}$ of
measurable functions : $\mathcal {C} (\mathbf{R}, \mathbf{R}_+)
\rightarrow \mathbf{R_+}$, such that : \\ \\
1) $F_0 = 0$, and $F_k \underset{k \rightarrow \infty}{\rightarrow}
F$. \\ \\
2) $(|F_k - F_{k-1} |) ((l^y)_{y \in \mathbf{R}})$ depends only on
$(l^y)_{|y| \leq c_k}$, and : $$ (|F_k - F_{k-1} |) ((l^y)_{y \in
  \mathbf{R}}) \leq \left( \frac{ \underset{|y| \leq c}{\sup} l^y +
    c_k} { \underset{|y| \leq c}{\inf} l^y + c_k} \right) ^n h_k
\left( \underset{|y| \leq c}{\inf} l^y \right)$$ 
\noindent 
3) $\underset{k \geq 1}{\sum} N_{c_k} (h_k) \leq M$. \\ \\ 
These conditions imply the following ones for the sequence $\left( G_k =
F_k^{(l_0^y)_{y \in \mathbf{R}}, x} \right)_{k \geq 1}$ : \\ \\
1) $G_0=  0$, and $G_k \underset{k \rightarrow \infty}{\rightarrow}
F^{(l_0^y)_{y \in \mathbf{R}}, x}$. \\ \\
2)  $(|G_k - G_{k-1} |) ((l^y)_{y \in \mathbf{R}})$ depends only on
$(l^z)_{|z| \leq c_k + |x|}$ and : 
$$(|G_k - G_{k-1} |) ((l^y)_{y \in \mathbf{R}}) \leq \left( \frac{
  \underset{z \in  [-c_k - x, c_k -x]}{\sup} (l_0^{z+x} + l^z) + c_k}{ 
 \underset{z \in  [-c_k - x, c_k -x]}{\inf } (l_0^{z+x} + l^z) + c_k}
  \right)^n h_k \left( \underset{z \in  [-c_k - x, c_k -x]}{\inf} 
(l_0^{z+x} + l^z) \right)$$ $$ \leq \left(  \frac{ \underset{z \in
  \mathbf{R}}{\sup} \, l_0^z +  \underset{|z| \leq c_k + |x|}{\sup} l^z +
  c_k + |x|} { \underset{|z| \leq c_k + |x|}{\inf}  l^z + c_k}
  \right)^n  h_k \left( \underset{ |z| \leq c_k + |x|   }{\inf} 
   l^z   \right)$$ $$ \leq 2^n \left[    \left(  \frac{  \underset{z \in
  \mathbf{R}}{\sup} \, l_0^z}{c_k} \right) ^n  + 
\left( \frac{  \underset{|z| \leq c_k + |x|}{\sup} l^z +
  c_k + |x|}{  \underset{|z| \leq c_k + |x|}{\inf} l^z +
  c_k + |x|} \right) ^n 
\left( 1 + \frac{|x|} { \underset{|z| \leq c_k + |x|}{\inf} l^z +
  c_k} \right) ^n \right]  h_k \left( \underset{ |z| \leq c_k + |x|   }{\inf} 
   l^z   \right)$$ $$ \leq 2^n \left(  \left( \underset{z \in
  \mathbf{R}}{\sup} \, l_0^z \right)^n + (1 + |x|)^n \right)
 \left( \frac{  \underset{|z| \leq c_k + |x|}{\sup} l^z +
  c_k + |x|}{  \underset{|z| \leq c_k + |x|}{\inf} l^z +
  c_k + |x|} \right) ^n  h_k \left( \underset{ |z| \leq c_k + |x|   }{\inf} 
   l^z   \right)$$
Therefore, $|G_k - G_{k-1}|$ satisfies the condition $C \left( c_k + |x|, n,
  2^n \left( \left( \underset{z \in \mathbf{R}}{\sup} \, l_0^z \right) ^n
  + (1 + |x|)^n \right) h_k \right)$. \\ \\
3) Now : $$ N_{c_k+|x|} \left(  2^n \left( \left( \underset{z \in
  \mathbf{R}}{\sup} \, l_0^z \right) ^n  + (1 + |x|)^n \right) h_k
  \right) $$ $$\leq  2^n \left( \left( \underset{z \in
  \mathbf{R}}{\sup} \, l_0^z \right) ^n  + (1 + |x|)^n \right)
  \frac{c_k+|x|}{c_k} N_{c_k} (h_k) $$ $$ \leq 2^n \left( 1 +  \left(
  \underset{z \in
  \mathbf{R}}{\sup} \, l_0^z \right) ^n \right) (1 + |x|)^{n+1}  N_{c_k}
  (h_k)$$
\noindent
and $\underset{k \geq 1}{\sum} N_{c_k}
  (h_k) \leq M$. \\ \\
Therefore : $$  N^{(n)} \left(  F^{(l_0^y)_{y \in \mathbf{R}}, x}
  \right)  \leq  2^n \left( 1 +  \left(
  \underset{z \in
  \mathbf{R}}{\sup} \, l_0^z \right) ^n \right) (1 + |x|)^{n+1} M$$ 
\noindent
By taking $M \rightarrow N^{(n)} (F)$, we obtain the majorization
stated in Lemma 5. \hfill$\Box$  \\ \\ 
\textbf{Proof of the Theorem : } $ \sqrt{2 \pi t} \, \mathbf{W} \left[ F
    \left( \left( l_t^y (X) \right)_{y \in \mathbf{R}} \right) \right]
    $ tends  to $I(F) > 0$ when $t$ goes to infinity, so it is
    strictly positive if $t$ is large enough, and $\mathbf{W}_t^F $ is
    well-defined. \\ \\
If $t$ is large enough, by Markov property : 
$$ \mathbf{W}_t^F (\Lambda_s) =  \mathbf{W} \left[
  \mathbf{1}_{\Lambda_s}  \frac{  \mathbf{W} \left[ F \left( \left( l_t^y
      (X) \right)_{y \in \mathbf{R}} \right) \left|  \sigma \{X_u, u \leq s
   \} \right. \right]} { \mathbf{W} \left[ F \left( \left( l_t^y
      (X) \right)_{y \in \mathbf{R}} \right) \right] } \right] $$ 
$$ =  \mathbf{W} \left[  \mathbf{1}_{\Lambda_s} 
 \frac{ \Psi_{t-s} \left( (l_s^y(X))_{y \in \mathbf{R}}, X_s \right) }  {\mathbf{W} \left[ F \left( \left( l_t^y
      (X) \right)_{y \in \mathbf{R}} \right) \right] } \right]$$
where, for all continuous functions $l$ from $\mathbf{R}$ to
 $\mathbf{R}_+$, and for all $x \in \mathbf{R}$, $u > 0$ : $$ \Psi_u \left(
 (l^y)_{y \in \mathbf{R}}, x \right) =  
 \mathbf{W} \left[  F^{(l^y )_{y \in
  \mathbf{R}}, x} \left( (l_{u}^y (X))_{y \in \mathbf{R}}
 \right) \right] $$
\noindent
By Proposition 4.2 : $$ \frac{  \Psi_{t-s} \left( (l_s^y(X))_{y \in
      \mathbf{R}}, X_s \right) }
 {\mathbf{W} \left[ F \left( \left( l_t^y
      (X) \right)_{y \in \mathbf{R}} \right) \right] } \underset{t
 \rightarrow \infty}{\rightarrow}  \frac{ I \left(  F^{(l_s^y (X))_{y \in
  \mathbf{R}}, X_s} \right)}{I(F)}$$
Moreover, for $t \geq 2s$ : $$ \sqrt{2 \pi t} \, \Psi_{t-s} \left(
  (l_s^y(X))_{y \in \mathbf{R}}, X_s \right) \leq \sqrt{\frac{t}{t-s}}
  N^{(n)}  \left(  F^{(l_s^y (X) )_{y \in
  \mathbf{R}}, X_s} \right) $$ $$\leq 
 2^{n+1/2} \left( 1 +
 \left( \underset{z \in \mathbf{R}}{\sup} \, l_s^z (X)  \right)^n \right)
 \left( 1 + |X_s| \right)^{n+1} N^{(n)} (F) $$ 
and for $t$ large enough : $$ \sqrt{2 \pi t} \, \mathbf{W} \left[ F
      \left( (l_t^y(X))_{y \in \mathbf{R}} \right) \right] \geq I(F)/2$$
\noindent
Hence, for $t$ large enough : $$ \frac{  \Psi_{t-s} \left( (l_s^y(X))_{y \in
      \mathbf{R}}, X_s \right) }
 {\mathbf{W} \left[ F \left( \left( l_t^y
      (X) \right)_{y \in \mathbf{R}} \right) \right] } \leq \frac{ 2^{n+3/2} \left( 1 +
 \left( \underset{z \in \mathbf{R}}{\sup} \, l_s^z (X)  \right)^n \right)
 \left( 1 + |X_s| \right)^{n+1} N^{(n)} (F)}{I(F)}$$
\noindent
Now : $$\mathbf{W} \left[  \left( 1 +
 \left( \underset{z \in \mathbf{R}}{\sup} \, l_s^z (X)  \right)^n \right)
 \left( 1 + |X_s| \right)^{n+1} \right] $$ $$\leq \left( \mathbf{W}
 \left[ \left( 1 +  
\left(  \underset{z \in \mathbf{R}}{\sup} \, l_s^z (X)  \right)^n
 \right)^2 \right] \right) ^{1/2} \left(  \mathbf{W} \left[  \left( 1 + |X_s|
 \right)^{2n+2} \right] \right)^{1/2} < \infty $$
since $  \underset{z \in \mathbf{R}}{\sup} \, l_s^z (X)$ and $|X_s|$ have
 moments of any order. \\ \\
By dominated convergence, we obtain the Theorem.  \hfill$\Box$
\section{Examples}
\noindent
In this section, we check that the conditions of the Theorem are
 satisfied in three examples studied by B. Roynette, P. Vallois and
 M. Yor, and one more particular case. \\ \\
I) \textbf{First example } (penalization with local time at level
 zero) \\ \\
 We take $F((l^y)_{y \in \mathbf{R}}) = \phi(l^0)$
 where $\phi$ is bounded and dominated by a positive, decreasing and
 integrable function $\psi$. \\ \\
$F$ satisfies the condition $C(1, 0, \psi)$. Hence : $$ 
N^{(0)} (F) \leq N_1 (\psi) = \psi(0) + \int_0^{\infty} \psi(y) dy  <
 \infty$$
On the other hand : $$ I(F) = 2 \int_0^{\infty} \phi(l) dl$$
$$ F^{(l_s^y (X))_{y \in
  \mathbf{R}}, X_s} ((l^y)_{y \in \mathbf{R}}) = l_s^0 (X) + l^{-X_s} $$
and : $$ I \left( F^{(l_s^y (X))_{y \in
  \mathbf{R}}, X_s}  \right) = \int_0^{\infty} dl \left( \mathbf{E}
  \left[ \phi (l_s^0 (X) + Y_{l,+}^{-X_s} ) \right]  + \mathbf{E}
  \left[ \phi (l_s^0 (X) + Y_{l,-}^{-X_s} ) \right]  \right) $$
Now, by using the fact that Lebesgue measure is invariant for BESQ(2)
  process, we obtain : 
$$  \int_0^{\infty} dl \, \mathbf{E}
  \left[ \phi (l_s^0 (X) + Y_{l,- sgn(X_s)}^{- X_s} ) \right]  =
  \int_0^{\infty} dl \, \phi (l_s^0 (X) + l) = \int_{l_s^0 (X)}^{\infty}
  \phi(l) dl$$
Moreover, the image of Lebesgue measure by a BESQ(0) process taken at
  time $x \geq 0$ is the sum of Lebesgue measure and $2x$ times Dirac
  measure at $0$;  more precisely, for all measurable functions $f :
  \mathbf{R}_+ \rightarrow \mathbf{R}_+$, one has : 
$$\int_0^{\infty} dl \, \mathbf{E} [f (Y_{l,-}^x)] = 2 x f(0) +
\int_0^{\infty} dy f(y) $$  Therefore : 
$$  \int_0^{\infty} dl \, \mathbf{E}
  \left[ \phi (l_s^0 (X) + Y_{l, sgn (X_s)}^{- X_s} ) \right] = 2
  |X_s| \phi (l_s^0 (X)) +  \int_{l_s^0 (X)}^{\infty}
  \phi(l) dl$$ and finally : 
$$I \left( F^{(l_s^y (X))_{y \in
  \mathbf{R}}, X_s}  \right) = 2 \left(   |X_s| \phi (l_s^0 (X)) 
+  \int_{l_s^0 (X)}^{\infty}
  \phi(l) dl \right) $$ 
Consequently, if $\phi$ is not a.e. equal to zero, we can apply 
the Theorem, and for $s \geq 0$, $\Lambda_s \in \mathcal{F}_s = \sigma
  \{X_u, u \leq s \}$ : $$  \mathbf{W}_{\infty}^F  (\Lambda_s) =
  \mathbf{W} \left(
  \mathbf{1}_{\Lambda_s} . \frac{|X_s| \phi( l_s^0(X)) + \Phi (l_s^0
  (X))}{\Phi (0)} \right) $$ where $\Phi(x) = \int_x^{\infty} \phi(l)
  dl $. \\ \\
This result is coherent with the limit measure obtained by
  B. Roynette, P. Vallois and M. Yor in \cite{7}. \\ \\ 
II) \textbf{Second example } (penalization with the supremum) \\ \\
We take  $F((l^y)_{y \in \mathbf{R}}) = \phi ( \inf \{y \geq 0, l^y =0 \}
  )$, where $\phi$ is dominated by a decreasing function
  $\psi : \mathbf{R}_+ \cup \{ \infty \} \rightarrow \mathbf{R}_+$
  such that $\int_0^{\infty} \psi (y) dy < \infty$. \\ \\
Let us recall that for this choice of $F$, $F \left( (l_t^y (X))_{y
  \in \mathbf{R}} \right)  = \phi(S_t)$, where $S_t$ denotes the
  supremum of $(X_s)_{s \leq t}$. \\ \\
Now, we take for $k \in \mathbf{N}$ : $$
F_k ((l^y)_{y \in \mathbf{R}}) = \phi_{2^k-1} (\inf \{ y \geq 0, l^y
  =0 \} )$$ where $ \phi_{2^k-1} = \phi. \mathbf{1}_{]- \infty,
  2^{k}-1[} $. \\ \\
1) One has $F_0 = 0$ and $F_k \underset{k \rightarrow
  \infty}{\rightarrow} F$ pointwise. \\ \\
2) $(|F_k - F_{k-1}|) ((l^y)_{y \in \mathbf{R}})$ depends only on
  $(l^y)_{|y| \leq 2^k - 1}$ and :  $$ (|F_k - F_{k-1}|) ((l^y)_{y \in
  \mathbf{R}}) \leq  \phi ( \inf \{y \geq 0, l^y =0 \} ) \mathbf{1}_{ 
\inf  \{y \geq 0, l^y =0 \} \in [2^{k-1} - 1, 2 ^k - 1[ } $$ $$ \leq \psi
  (2^{k-1} - 1) \mathbf{1}_{ \underset{|y| \leq 2^k-1}{\inf} l^y =
  0}$$
Hence, $|F_k - F_{k-1}|$ satisfies the condition $C(2^k-1, 0,
  \psi(2^{k-1} -1) \mathbf{1}_{ \{ 0 \} })$. \\ \\
3) Therefore : $$ N^{(0)} (F) \leq \underset{k \geq 1}{\sum}  (2^k-1)
  \psi (2^{k-1} -1) \leq \psi(0) + 4 \int_0^{\infty} \psi < \infty$$
Moreover : 
$$I(F) = \int_0^{\infty} dl \,\mathbf{E} \left[ \phi \left( \inf \{ y
    \geq 0, Y_{l,+}^y = 0 \} \right) \right] + 
 \int_0^{\infty} dl \, \mathbf{E} \left[ \phi \left( \inf \{ y
    \geq 0, Y_{l,-}^y = 0 \} \right) \right]$$
The first integral is equal to zero and $\inf \{ y
    \geq 0, Y_{l,-}^y = 0 \}$ is the inverse of an exponential
    variable of parameter $l/2$. \\ 
Therefore : 
$$ I(F) = \int_0^{\infty} dl \int_0^{\infty} dy \, \frac{l}{2 y^2}
e^{-l/2y} \phi(y) dy $$ $$= \int_0^{\infty} dy \, \phi(y) \int_0^{\infty} dl
\, \frac{l}{2 y^2} e^{-l/2y} = 2 \int_0^{\infty} \phi(y) dy $$
By similar computations, we obtain : 
$$ I \left(  F^{(l_s^y (X))_{y \in
  \mathbf{R}}, X_s} \right) = \int_0^{\infty} dl \, \mathbf{E} [ \phi
  (S_s \vee (X_s + \inf \{ y
    \geq 0, Y_{l,-}^y = 0 \} ))] $$ $$ = 2 \left( (S_s - X_s) \phi(S_s) +
  \int_{S_s}^{\infty} \phi(y) dy \right) $$
Consequently, if $\phi$ is not a.e. equal to zero, the sequence
  $(\mathbf{W}_t^F)_{t \geq 0}$ satisfies for every $s \geq 0$,
  $\Lambda_s \in \mathcal{F}_s =  \sigma \{ X_u, u \leq s \} $ : 
$$ \mathbf{W}_t^F (\Lambda_s) \underset{t \rightarrow
  \infty}{\rightarrow} \mathbf{W}_{\infty}^F (\Lambda_s)$$ 
where  $$ \mathbf{W}_t^F = \frac{\phi(S_t)}{\mathbf{W} [ \phi(S_t)]
  } . \mathbf{W} $$ 
and : $$ \mathbf{W}_{\infty}^F (\Lambda_s) = \mathbf{W} \left[
  \mathbf{1}_{\Lambda_s} \frac{ (S_s - X_s) \phi(S_s) + \Phi(S_s)}
  {\Phi (0)}  \right] $$
It corresponds to B. Roynette, P. Vallois and M. Yor's penalization
  results for the supremum (see \cite{7}). \\ \\
III) \textbf{Third example } (exponential penalization with an integral
  of the local times) \\ \\
Let us take : $F((l^y)_{y \in \mathbf{R}}) = \exp \left( - \int_{-\infty}^
 {\infty} V(y) l^y dy \right)$, where $V$ is a positive measurable
  function, not a.e. equal to zero, and integrable with respect to
 $(1 + y^2) dy$ (this condition is a little more restrictive than the
  condition obtained by B. Roynette, P. Vallois and M. Yor in \cite{6}). \\ \\
In that case, there exists $c \geq 1$ such that : $$ \int_{-c}^c
 V(y) dy > 0$$ We consider the following approximations of $F$ : \\ 
$F_0 = 0$, and for $k \geq 1$, $F_k  ((l^y)_{y \in \mathbf{R}}) = \exp
  \left( - \int_{-2^k c}^{2^k c} V(y) l^y dy \right)$. \\ 
The following holds : \\ \\
1) $F_0 = 0$ and $F_k \underset{k \rightarrow \infty}{\rightarrow}
  F$. \\ \\
2) $|F_k - F_{k-1}| ((l^y)_{y \in \mathbf{R}})$ depends only on
  $(l^y)_{y \in [-2^k c, 2^k c]}$ and if $k \geq 2$ : 
$$ |F_k - F_{k-1}| ((l^y)_{y \in \mathbf{R}}) \leq \left( \underset{[-2^k c,
  2^k c] \backslash [-2^{k-1} c, 2^{k-1} c]}{\int} V(y) dy
  \right)...$$ $$...
  \left( \underset{y \in [-2^k c, 2^k c]}{\sup} l^y \right) \exp
  \left(   - \int_{-2^{k-1} c}^{2^{k-1} c} V(y) l^y dy \right)  $$
$$ \leq \left( \underset{[-2^k c,
  2^k c] \backslash [-2^{k-1} c, 2^{k-1} c]}{\int} V(y) dy \right)
 \left( \frac{  \underset{y \in [-2^k c, 2^k c]}{\sup} l^y + 2^k c} 
{  \underset{y \in [-2^k c, 2^k c]}{\inf} l^y + 2^k c} \right)
\left( \underset{y \in [-2^k c, 2^k c]}{\inf} l^y + 2^k c \right) 
 ... $$ $$... \exp
  \left[   - \left( \int_{-2^{k-1} c}^{2^{k-1} c} V(y) dy \right) 
 \left(  \underset{y \in [-2^k c, 2^k c]}{\inf} l^y \right) \right]$$ 
Moreover :
$$ |F_1 - F_0|  ((l^y)_{y \in \mathbf{R}}) \leq \exp \left[ - \left(
    \int_{2c}^{2c} V(y) dy \right) \left(  \underset{y \in [-2 c,
    2 c]}{\inf} l^y \right) \right] $$
Therefore, if we put $ \rho = \int_{-c}^c V(y) dy > 0$, 
 for every $k \geq 1$, $|F_k - F_{k-1}|$ satisfies the
    condition $C(2^k c, 1, h_k)$ where the decreasing function $h_k$
    is defined by : 
$$ h_k (l) =  \left(\mathbf{1}_{k = 1} +  \underset{[-2^k c,
  2^k c] \backslash [-2^{k-1} c, 2^{k-1} c]}{\int} V(y) dy \right)
  (l + 2^k c + \rho^{-1}) e^{- \rho l}$$ 
\noindent
3) One has : 
$$N_{2^k c} (h_k) \leq   \left(\mathbf{1}_{k = 1} +  \underset{[-2^k c,
  2^k c] \backslash [-2^{k-1} c, 2^{k-1} c]}{\int} V(y) dy \right)
 ( 2^{2k} c^2 + 2^{k+1} c \rho^{-1} + 2 \rho^{-2})$$
Hence : 
 $$\underset{k \geq 1}{\sum} N^{(1)}_{2^k c} (h_k) \leq 
(1 + \rho^{-1} + \rho^{-2}) \left( 4c^2 +  \underset{k \geq 1}{\sum} 
2^{2k} c^2 \underset{[-2^k c,
  2^k c] \backslash [-2^{k-1} c, 2^{k-1} c]}{\int} V(y) dy \right) $$
$$ \leq 4 (1 + \rho^{-1} + \rho^{-2})  \left( c^2 + \int_{\mathbf{R}} (1 +
  y^2) V(y) \right)<\infty $$
Moreover, by properties of BESQ processes, for all $l \geq 0$, $y \in
\mathbf{R}$ : $$\mathbf{E} \left[ Y_{l,+}^y \right] \leq l + 2 |y|$$
and $$ \mathbf{E} \left[ \int_{\mathbf{R}}  Y_{l,+}^y V(y) dy
\right] \leq \int_{\mathbf{R}} (l+2|y|) V(y) dy < \infty$$
Therefore : $$\mathbf{E} \left[ \exp  \left( -  \int_{\mathbf{R}}
    Y_{l,+}^y V(y) dy \right) \right] > 0$$ and $I(F) > 0$. \\ \\
Consequently, the Theorem applies in this case and B. Roynette,
  P. Vallois and M. Yor's penalization result holds (see \cite{6}).\\ \\
IV) \textbf{Fourth example } (penalization with local times at two
  levels) \\ \\  This example is a generalization of the
  first one. \\ Let us take, for $y_1 < y_2$,  $F((l^y)_{y \in \mathbf{R}}) =
  \phi(l^{y_1}, l^{y_2})$  where $\phi (l_1, l_2) \leq h (l_1 \wedge
  l_2)$ for a positive, integrable and decreasing function $h$. \\ \\
In that case, $F$ satisfies the condition $C(|y_1| \vee |y_2|, 0, h)$,
  so the Theorem applies if we have $I(F) > 0$. \\ \\
For $y>0$, $z$, $z' \geq 0$, let $q^{(0)}_y (z, z')$ be the density at
  $z'$ of a BESQ(0) process starting from level $z$ and taken at time
  $y$, $Q^{(0)}_y (z,0)$ the probability that this process is equal to
  zero, and $q^{(2)}_y (z, z')$ the density at
  $z'$ of a BESQ(2) process starting from $z$ and taken at time
  $y$. If $0 < y_1 < y_2$, one has : $$ I(F) = \int_0^{\infty} dl
  \int_0^{\infty} dl_1  \int_0^{\infty} dl_2 \, q^{(2)}_{y_1} (l, l_1)
  \,  q^{(2)}_{y_2-y_1} (l_1, l_2) \, \phi (l_1, l_2)$$ $$ +  \int_0^{\infty} dl
  \int_0^{\infty} dl_1   \int_0^{\infty} dl_2 \, q^{(0)}_{y_1} (l, l_1)
  \,  q^{(0)}_{y_2-y_1} (l_1, l_2) \, \phi (l_1, l_2) $$ $$+
  \int_0^{\infty} dl \int_0^{\infty} dl_1  \, q^{(0)}_{y_1} (l, l_1)
  \,  Q^{(0)}_{y_2-y_1} (l_1, 0) \, \phi(l_1, 0) +  \int_0^{\infty} dl
  \,  Q^{(0)}_{y_1} (l, 0) \phi(0,0) $$
\noindent
Now, by properties of time-reversed BESQ processes : $q^{(0)}_y (z, z') =
 q^{(4)}_y (z', z)$ (where $q^{(4)}$ is the density of the BESQ(4)
 process) and $q^{(2)}_y (z, z') = q^{(2)}_y (z', z)$. Hence : 
$$ \int_0^{\infty} q^{(0)}_y (z, z') \, dz =  \int_0^{\infty} q^{(4)}_y
(z', z) \, dz =  1$$ and $$ \int_0^{\infty} q^{(2)}_y (z, z') \, dz =
\int_0^{\infty} q^{(2)}_y
(z', z) \, dz =  1$$ since $q^{(2)}$ and $q^{(4)}$ are probability
densities with respect to the second variable. \\ \\
Moreover : $$ \int_0^{\infty} Q_y^{(0)} (z,0) dz =   \int_0^{\infty}
e^{-z/2y} dz  = 2y$$ Therefore : $$ I(F) =  \int_0^{\infty} dl_1
\int_0^{\infty} dl_2 \, ( q^{(2)}_{y_2-y_1} (l_1,l_2) +
q^{(0)}_{y_2-y_1} (l_1,l_2)) \phi (l_1, l_2) $$ $$ +  \int_0^{\infty}
dl_1 \,  Q^{(0)}_{y_2-y_1} (l_1, 0) \phi(l_1,0) + 2 y_1 \phi(0,0) $$ 
\noindent
for $0 \leq y_1 < y_2$. \\ \\
Similar computations give for $y_1  < y_2 \leq 0$ : 
$$ I(F) =  \int_0^{\infty} dl_1
\int_0^{\infty} dl_2 \, ( q^{(2)}_{y_2-y_1} (l_2,l_1) +
q^{(0)}_{y_2-y_1} (l_2,l_1)) \phi (l_1, l_2) $$ $$ +  \int_0^{\infty}
dl_2 \,  Q^{(0)}_{y_2-y_1} (l_2, 0) \phi(0, l_2) + 2 |y_2| \phi(0,0) $$ 
\noindent
For $y_1 < 0 < y_2$, we have : $$ I(F) =\int_0^{\infty} dl
\int_0^{\infty} dl_1  \int_0^{\infty} dl_2 \,  q^{(2)}_{y_2} (l, l_2)
  \,  q^{(0)}_{|y_1|} (l, l_1) \, \phi (l_1, l_2) $$ $$+ 
\int_0^{\infty} dl \int_0^{\infty} dl_2 \,  q^{(2)}_{y_2} (l, l_2)
  \,  Q^{(0)}_{|y_1|} (l, 0) \, \phi (0, l_2)$$ $$ + 
\int_0^{\infty} dl
\int_0^{\infty} dl_1  \int_0^{\infty} dl_2 \,  q^{(0)}_{y_2} (l, l_2)
  \,  q^{(2)}_{|y_1|} (l, l_1) \, \phi (l_1, l_2) $$ $$ + 
\int_0^{\infty} dl \int_0^{\infty} dl_1 \,  Q^{(0)}_{y_2} (l, 0)
  \,  q^{(2)}_{|y_1|} (l, l_1) \, \phi (l_1, 0)$$
\noindent
Now, for $y'$, $y''>0$, and $z$, $z'$, $z'' \leq 0$, the two following
equalities hold : 
$$\int_0^{\infty}  q^{(2)}_{y'} (z, z')   q^{(0)}_{y''} (z, z'') \, dz
= \frac{ y' q^{(2)}_{y'+y''} (z', z'') + y'' q^{(0)}_{y' + y''}
  (z',z'')} {y' + y''}$$ 
$$ \int_0^{\infty}  q^{(2)}_{y'} (z, z')   Q^{(0)}_{y''} (z, 0) \, dz
= \frac{  y''} {y' + y''}  Q^{(0)}_{y' + y''}
  (z',0) $$ 
\noindent
(the first one can be proven by using \cite{9}, Lemma 3, and the relation
: $q^{(0)}_y (z,z') = q^{(4)}_y (z',z)$; the second is a consequence
of the equality : $Q^{(0)}_{y''} (z,0) = e^{-z/2y''} = 2
y'' q^{(2)}_{y''} (0,z)$).     \\
\\
Therefore : 
$$ I(F) = \int_0^{\infty} dl_1 \int_0^{\infty} dl_2 \left[
  q^{(2)}_{y_2-y_1} (l_1,l_2)  + \frac{ |y_1|  q^{(0)}_{y_2-y_1}
  (l_2,l_1) + y_2  q^{(0)}_{y_2-y_1}
  (l_1,l_2)}{y_2 - y_1} \right] \phi(l_1,l_2)  $$
$$+  \int_0^{\infty} dl_1 \, \frac{y_2}{y_2 - y_1} \, Q^{(0)}_{y_2 -
  y_1} (l_1,0) \, \phi(l_1,0)$$ $$ +  \int_0^{\infty} dl_2 \,
  \frac{|y_1|}{y_2 - y_1} \, Q^{(0)}_{y_2 -
  y_1} (l_2,0) \, \phi(0,l_2)$$
\noindent
This computation of $I(F)$ implies the following : \\ \\
1) For $0<y_1<y_2$, the Theorem applies iff : 
 $$ \int_0^{\infty} dl_1 \int_0^{\infty} dl_2 \, \phi(l_1,l_2) +
 \int_0^{\infty} dl_1  \, \phi(l_1,0)  + \phi(0,0) > 0$$
\noindent
2) For $0=y_1<y_2$, it applies iff :
 $$ \int_0^{\infty} dl_1 \int_0^{\infty} dl_2 \, \phi(l_1,l_2)  +
 \int_0^{\infty} dl_1  \, \phi(l_1,0) > 0$$
\noindent
3) For $y_1 < 0 < y_2$, it applies iff : 
 $$ \int_0^{\infty} dl_1 \int_0^{\infty} dl_2 \, \phi(l_1,l_2)
 +  \int_0^{\infty} dl_1  \, \phi(l_1,0)  +
 \int_0^{\infty} dl_2  \, \phi(0,l_2) > 0$$ 
\noindent
4) For $y_1 < y_2 = 0$, it applies iff :
 $$ \int_0^{\infty} dl_1 \int_0^{\infty} dl_2 \, \phi(l_1,l_2) + 
 \int_0^{\infty} dl_2  \, \phi(0,l_2) > 0$$ 
\noindent
5) For $y_1 < y_2 < 0$, it applies iff : 
 $$ \int_0^{\infty} dl_1 \int_0^{\infty} dl_2 \, \phi(l_1,l_2)+ 
\int_0^{\infty} dl_2  \, \phi(0,l_2) + \phi(0,0) > 0$$
\noindent
If the Theorem holds, it is possible to compute $I(F^{(l_s^y(X))_{y
    \in \mathbf{R}}, X_s})$ in order to obtain the density, restricted
    to $\mathcal{F}_s$,  of
    $\mathbf{W}_{\infty}^{F}$ with respect to $\mathbf{W}$. \\ \\
For $X_s \leq y_1 < y_2$, we have : 
 $$I(F^{(l_s^y(X))_{y
    \in \mathbf{R}}, X_s})  =  \int_0^{\infty} dl_1
\int_0^{\infty} dl_2 \, ( q^{(2)}_{y_2-y_1} (l_1,l_2) +
q^{(0)}_{y_2-y_1} (l_1,l_2)) \phi (l_s^{y_1} (X) + l_1,l_s^{y_2} (X) + l_2) $$ $$ +  \int_0^{\infty}
dl_1 \,  Q^{(0)}_{y_2-y_1} (l_1, 0) \phi(l_s^{y_1} (X) + l_1,l_s^{y_2} (X)  ) +
 2 (y_1-X_s) \phi(l_s^{y_1} (X),l_s^{y_2} (X)    ) $$ 
\noindent
For $y_1 < X_s < y_2$ : 
$$I(F^{(l_s^y(X))_{y
    \in \mathbf{R}}, X_s})  = \int_0^{\infty} dl_1 \int_0^{\infty} dl_2 \left[
  q^{(2)}_{y_2-y_1} (  l_1, l_2    ) ...\right. $$ $$ \left. ... 
  + \frac{ (X_s - y_1)  q^{(0)}_{y_2-y_1}
  (l_2,l_1) + (y_2-X_s)   q^{(0)}_{y_2-y_1}
  (l_1,l_2)}{y_2 - y_1} \right] \phi( l_s^{y_1} (X) + l_1,l_s^{y_2}
    (X) + l_2       )  $$
$$+  \int_0^{\infty} dl_1 \, \frac{y_2-X_s}{y_2 - y_1} \, Q^{(0)}_{y_2 -
  y_1} (l_1,0) \, \phi(l_s^{y_1} (X) + l_1, l_s^{y_2} (X))$$ $$ +  \int_0^{\infty} dl_2 \,
  \frac{X_s-y_1}{y_2 - y_1} \, Q^{(0)}_{y_2 -
  y_1} (l_2,0) \, \phi( l_s^{y_1} (X),l_s^{y_2} (X) + l_2         )$$ 
\noindent
For $y_1<y_2 \leq X_s$ : 
 $$I(F^{(l_s^y(X))_{y
    \in \mathbf{R}}, X_s})  =  \int_0^{\infty} dl_1
\int_0^{\infty} dl_2 \, ( q^{(2)}_{y_2-y_1} (l_2,l_1) +
q^{(0)}_{y_2-y_1} (l_2,l_1)) \phi (l_s^{y_1} (X) + l_1,l_s^{y_2} (X) + l_2) $$ $$ +  \int_0^{\infty}
dl_2 \,  Q^{(0)}_{y_2-y_1} (l_2, 0) \phi(l_s^{y_1} (X),l_s^{y_2} (X) +
    l_2) +
 2 (X_s - y_2) \phi(l_s^{y_1} (X),l_s^{y_2} (X)    ) $$
\noindent
These formulae give an explicit expression for the limit measure
obtained in our last example. \\ \\
\textbf{Remark 5.1 : } It is not difficult to extend this example to a
functional of a finite number of local times. We have only considered
the case of two local times in order to avoid too complicated
notation. \\ \\
\textbf{Remark 5.2 : } The main Theorem cannot be extended to every
functional $F$. For example, if we consider the functional :
$$F((l^y)_{y \in \mathbf{R}}) = \exp \left( - \int_{- \infty}^{\infty}
  (l^y)^2 \, dy \right)$$ 
\noindent
which corresponds to Edwards' model in dimension 1 (see \cite{10}),  
the expectation $ \mathbf{E} [ F((L_t^y)_{y \in \mathbf{R}})]$ tends
exponentially to zero, and $I(F) = 0$, since for all $l>0$ : $$\int_{-
  \infty}^{\infty} (Y_{l,+}^y)^2 \, dy = \infty$$ almost surely. \\ \\
Therefore, it is impossible to study this case as the examples given
above. \\ 
Another case for which the Theorem cannot apply is the functional : 
$$F((l^y)_{y \in \mathbf{R}}) = \phi ( \sup (l^y)_{y \in \mathbf{R}}
)$$ where $\phi$ is a bounded function with compact support. \\ \\
It would be interesting to find another way to study this kind of
penalizations. 
\bibliographystyle{alpha}
 \bibliography{locaux}

\end{document}